\documentclass[12pt]{article}
\usepackage{amsthm}
\usepackage{amssymb}
\usepackage{graphicx}
\usepackage{float}
\usepackage{enumerate}
\usepackage{setspace}
\usepackage{perpage}
\usepackage[compact]{titlesec}
\usepackage{longtable}
\usepackage{mathtools}
\usepackage{nomencl}
\makeglossary

\theoremstyle{plain}

\begin{document}

 \baselineskip 0.22in

 \title{COMPARISON OF SEVERAL REWEIGHTED $l_1$-ALGORITHMS FOR SOLVING CARDINALITY MINIMIZATION PROBLEMS}
\author {  Mohammad Javad Abdi \thanks{School of Mathematics, University of Birmingham, Edgbaston B15 2TT, United Kingdom. Email: abdimj@maths.bham.ac.uk}  }

%\date{}

\maketitle

\textbf{Abstract.}
Reweighted $l_1$-algorithms have attracted a lot of attention in the field of applied mathematics. A unified framework of such algorithms has been recently proposed in \cite{1}. In this paper we construct a few new examples of reweighted $l_1$-methods. These functions are certain concave approximations of the $l_0$-norm function. We focus on the numerical comparison between some new and existing reweighted $l_1$-algorithms. We show how the change of parameters in reweighted algorithms may affect the performance of the algorithms for finding the solution of the cardinality minimization problem. In our experiments, the problem data were generated according to different statistical distributions, and we test the algorithms on different sparsity level of the solution of the problem. Our numerical results demonstrate that the reweighted $l_1$-method is one of the efficient methods for locating the solution of the cardinality minimization problem.

\newpage
\section{Introduction}

The cardinality minimization problem over a convex set $C\subset \mathbb{R}^n$ can be written as
\begin{eqnarray}\begin{split}\label{CMPMAIN}
\mbox{Minimize } &\|x\|_0& \\
\mbox {s.t.  } & x\in C.&\\
\end{split}\end{eqnarray}
This problem is to minimize the number of non-zero components of a vector satisfying certain constraints. In other words, cardinality minimization problem is looking for the sparsest vector in a given feasible set.

In this paper, we suppose that $C$ is defined by an undetermined system of linear equations, i.e,
$$C=\{x; Ax=b\},~\mbox{where}~A\in \mathbb{R}^{m\times n}~(m<n)~,b\in \mathbb{R}^{m}.$$
These linear systems have infinite many solutions, and the purpose of cardinality minimization problem(CMP) is to find the sparsest one, which can be stated as follows:
\begin{eqnarray}\begin{split}\label{CMP}
\mbox{Minimize } &\|x\|_0& \\
\mbox {s.t.  } & Ax=b.&\\
\end{split}\end{eqnarray}

The problem ($\ref{CMP}$) is closely related to compressed sensing which is dealing with the reconstruction of sparse signals from a limited number of linear measurements \cite{34,35,36,37}. Also problems with cardinality constraints have a wide range of applications, especially in portfolio optimization problems \cite{39,40}, and principal component analysis and model reduction \cite{41,43}. The more generalized version of cardinality minimization problems is so called rank minimization problems which have been considered in recent years \cite{8,9}.

The $l_0$-norm function is discontinuous, so the main idea for solving the problem ($\ref{CMP}$) is to approximate the $l_0$-norm function by some other continuous functions which are easier to deal with. For example $l_p$-norm function($0<p<1$) is one of the approximations of the $l_0$-norm. $l_p$ minimization $(0<p<1)$ has been studied in \cite{44,46,26}. Figure (1) represents the graph of $\sum_{i=1}^n|x_i|^p$, for $p=1$, $p=0.6$, and $p=0.2$. Note that as $p$ goes to zero, $\sum_{i=1}^n|x_i|^p$ approaches to the $l_0$-norm function.
\begin{figure}[htp]
\centering
$\begin{array}{c}
\includegraphics[width=0.5\textwidth,totalheight=0.25\textheight]{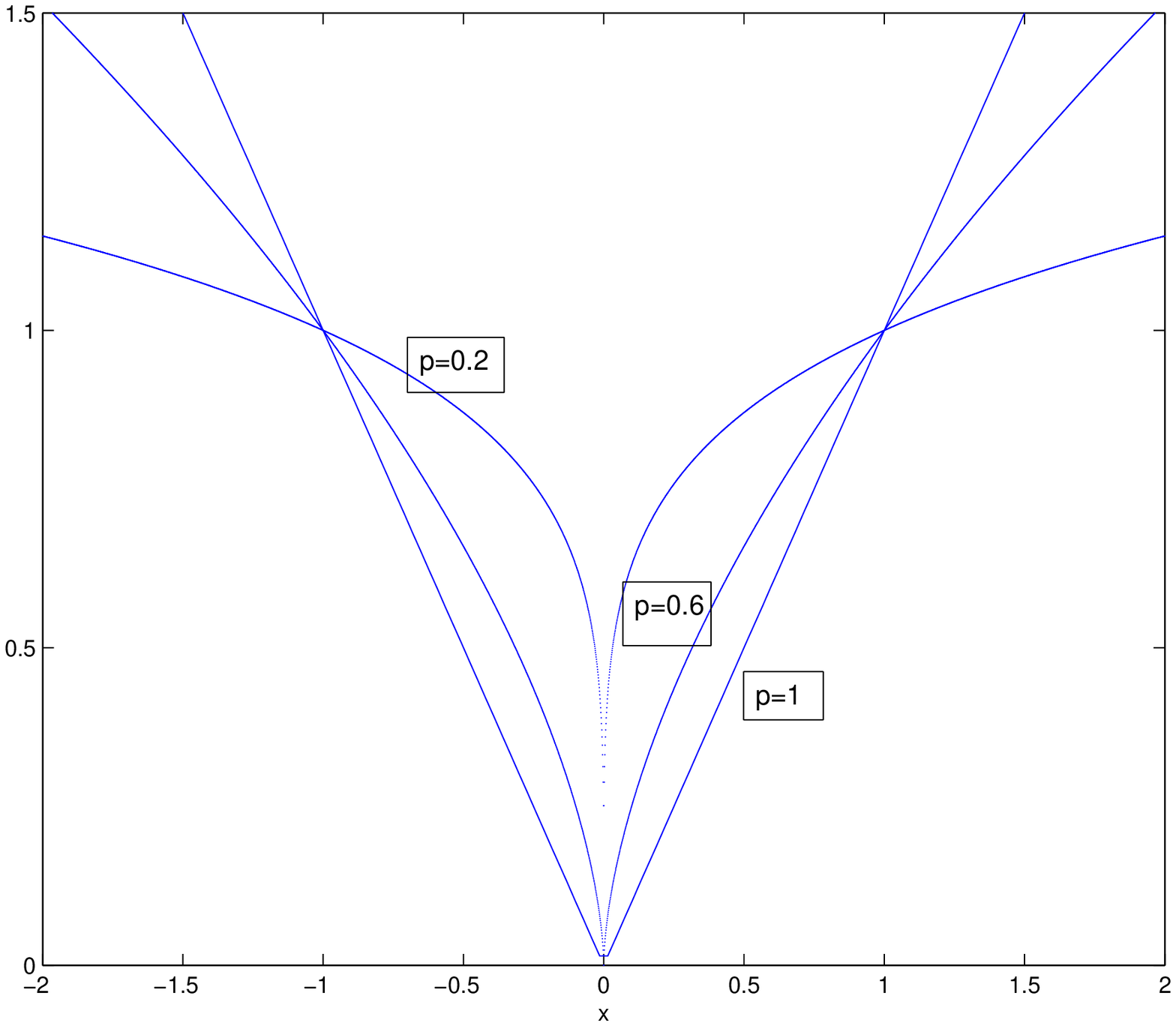}
\label{FIG1}
\end{array}$
\caption{\textrm{The graph of} $\sum_{i=1}^n|x_i|^p$ \textrm{for different values of} $0<p<1.$ }
$\begin{array}{c}
\includegraphics[width=0.5\textwidth,totalheight=0.25\textheight]{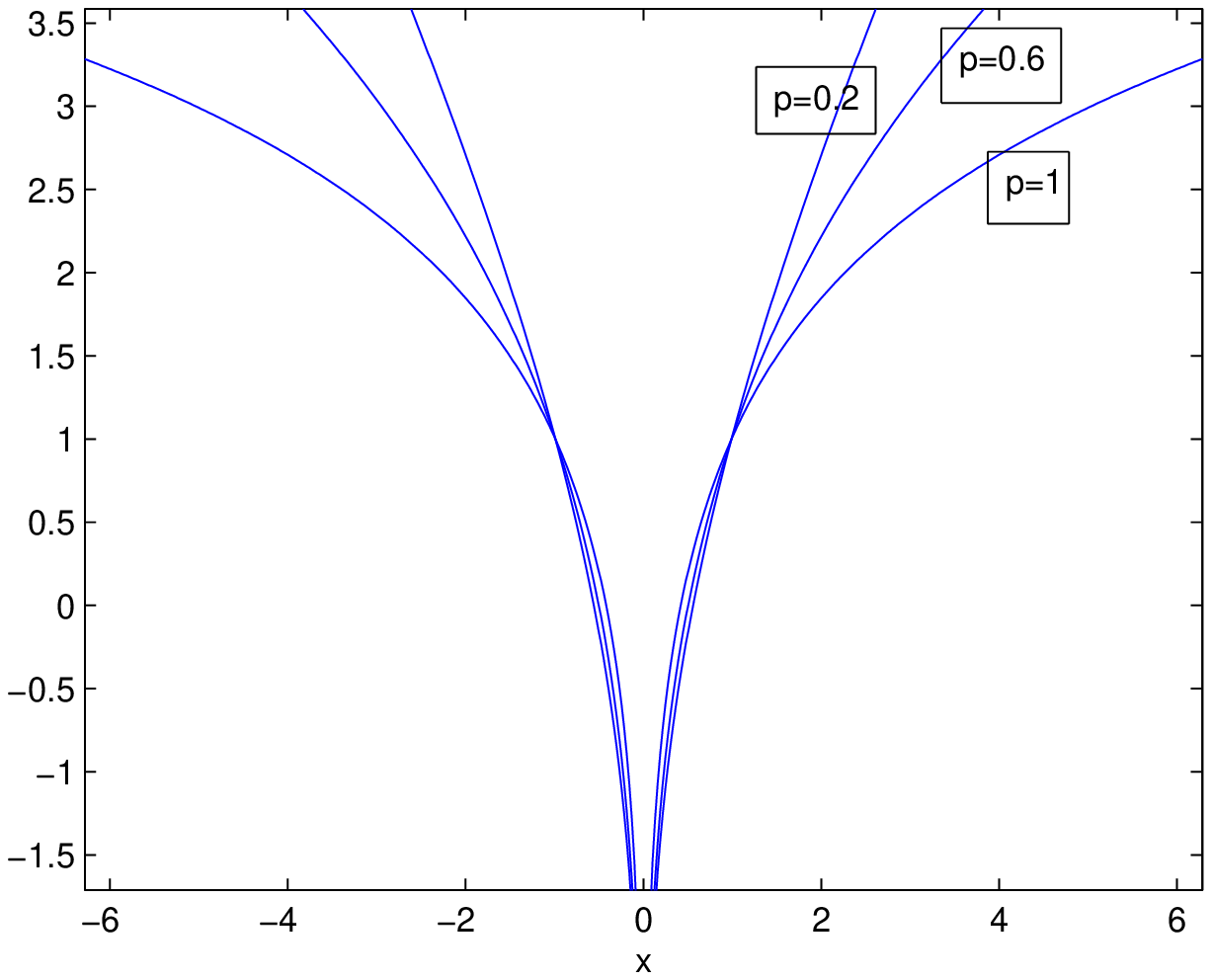}
\label{FIG2}
\end{array}$
\caption{\textrm{The graph of} $\sum_{i=1}^n \log(|x_i|+\epsilon)+\sum_{i=1}^n(|x_i|+\epsilon)^p$ \textrm{for different values of} $0<p<1.$ }
\end{figure}
As seen in the Figure (1), the closest convex approximation of $\|x\|_0$ is the well known $l_1$-norm function. So it is
unavoidable to use non-convex functions, especially concave functions, in order to have a better approximation of $\|x\|_0$. In \cite{1} Zhao and Li introduced the following function which is a combination of $l_p$-norm$(0<p<1)$ and the $\log$ function to approximate $\|x\|_0$,(see Figure (2)):
$$F_{\epsilon}(x)=\sum_{i=1}^n \log(|x_i|+\epsilon)+\sum_{i=1}^n(|x_i|+\epsilon)^p,~0<p<1.$$
We will discuss this kind of functions later. Finite successive linear approximation algorithms have also been used and are still used to get an approximated solution of the concave approximation problems \cite{18,19,20,60}.

$l_1$-minimization methods have been used to solve the problem ($\ref{CMP}$). This is motivated by the main idea of replacing $\|x\|_0$ function with its local convex envelop, the $l_1$-norm function, and then solve the resulting linear program \cite{3,4,48}. Under certain conditions the $l_1$-minimization method is able to obtain the exact solution of the problem  ($\ref{CMP}$) for very sparse solutions of the system $Ax=b$. In the literature, several conditions have been introduced and discussed for the equivalence between the $l_1$-minimization and the $l_0$-minimization. The outstanding ones are spark \cite{4}, mutual coherence \cite{5,6}, restricted isometry property(RIP)\cite{3,10,11}, and null space property(NSP) \cite{12,13}. For large optimization problems, the unconstrained version of the problem has been investigated in the literature, that may be referred as Lasso-type problems \cite{52}.

Numerical experiments show that weighted approaches are very affective in locating an exact solution of the problem ($\ref{CMP}$), and it can outperform other methods, in many situations \cite{1,22,23,24,25,26}.\\

Candes, Wakin, and Boyd \cite{23} proposed a weighted $l_1$-algorithm as follows:
\begin{eqnarray}\begin{split}\label{WP1}
\mbox{Minimize } &\sum_{i=1}^n\omega_i|x_i|& \\
\mbox {s.t.  } & Ax=b.&\\
\end{split}\end{eqnarray}
By introducing a diagonal matrix $W=\mbox{diag}(\omega_1,\omega_2,...,\omega_n)$, the problem above can be written as
\begin{eqnarray}\begin{split}\label{WP2}
\mbox{Minimize } &\|Wx\|_1& \\
\mbox {s.t.  } & Ax=b.&\\
\end{split}\end{eqnarray}
The weight can be interpreted as penalties for the components of the vector $x$. Larger penalties, ($\omega_i$)s, apply to smaller component of the vector $x$, for example one may choose the weights as $\frac{1}{|x_i|},~i=1,...,n$.
However to avoid having infinity penalties, one may add a parameter like $\epsilon>0$ to define the following weight \cite{23}
$$\omega_i=\frac{1}{|x_i|+\epsilon},~i=1,...,n.$$

Choosing a proper $\epsilon$ to have a more efficient algorithm is one of the challenges. Very small or very large $\epsilon$ might lead to improper weights which may cause the failure of the algorithms. Since very small $\epsilon$ might result in infinity penalties for small component of $x$, and with very big $\epsilon$ the penalty might not recognize the difference between the small components of $x$ and the large ones. We discuss the choice of $\epsilon$ for our algorithms in the numerical experiment later.\\

In an iterative reweighted $l_1$-algorithm, weights can be defined from the iteration in the previous step. Suppose the solution at the step $l$ is $x^l$, then the weight at the next step $l+1$, can be given as $\omega^{l+1}=\frac{1}{|x^l|+\epsilon}$. This was introduced by Candes, Wakin, and Boyd, and we refer to their algorithm as CWB, in this paper.\\

In \cite{1}, Zhao and Li introduced a unified framework for the reweighted $l_1$-minimization. The main idea is to define a merit function which is a certain concave approximation of the cardinality function, and to construct different types of weights through the linearization techniques. Based on the class of merit functions defined by Zhao and Li \cite{1}, we identify several new specific merit functions which are used to define the weights of the reweighted $l_1$-algorithms. The main purpose of this paper is to study these merit functions, and to test the success probability of the reweighted algorithms associated with these merit functions for locating the sparsest solution of linear systems, where the matrices, A, are generated based on different statistical distributions. Note that most of the previous experiments in the literature use normally distributed matrices. Also, we demonstrate how the parameters used in the algorithms may affect the performance of these methods. Furthermore, we evaluate what choice of $\epsilon$ may make our algorithms work better. In section 2, we discuss different types of merit functions and the associated reweighted $l_1$-algorithms. In section 3, we present and discuss our numerical results, and provide comparison between these algorithms.
\section{Merit functions and reweighted algorithms}
Merit functions have been used frequently in the field of optimization. Recently Zhao and Li \cite{1} has used merit functions to approximate $l_0$-norm. The merit function is defined as follows.

\textbf{Merit function:} For any $\epsilon>0$, a merit function $F_{\epsilon}(x):\mathbb{R}^n\rightarrow \mathbb{R}$ for approximating the $l_0$-norm, is strictly concave, separable, coercive, strictly increasing and twice differentiable, with the following properties:
\begin{enumerate}
  \item $\lim_{\epsilon\rightarrow 0}\frac{F_{\epsilon}(x)}{g(\epsilon)}=\|x\|_0+C,~g(\epsilon)>0~ \mbox{is a function of }\epsilon, ~\mbox{and}~C~\mbox{is a constant},$\\

  \item $F_{\epsilon}(x)=F_{\epsilon}(|x|),~~\forall x\in \mathbb{R}^n,$\\

  \item $\lim_{(x_i,\epsilon)\rightarrow (0,0)}[\nabla F_{\epsilon}(x)]_i=\infty,~~\forall x\geq 0,~\forall\epsilon>0, $\\

  \item $\lim_{\epsilon \rightarrow 0}[\nabla F_{\epsilon}(x)]_i=c_i,~~\forall x_i>0,~\mbox{where each}~c_i~\mbox{is a positive constant}.$
\end{enumerate}

After replacing the $\|x\|_0$ by a merit function the problem ($\ref{CMP}$) can be written as
\begin{eqnarray}\begin{split}\label{Merit}
\mbox{Minimize } &F_{\epsilon}(x)& \\
\mbox {s.t.  } & Ax=b.&\\
\end{split}\end{eqnarray}

Note that $F_{\epsilon}(x)$ is a concave function. One of the usual methods to solve concave optimization problems is to apply the linearization method, which in this case is a special type of Majorization-Minimization(MM) method. For more illustration see that by applying the Taylor expansion of $F_{\epsilon}(x)$ around a point $u$, we conclude
$$F_{\epsilon}(x)\leq F_{\epsilon}(u)+\langle\nabla F_{\epsilon}(u),x-u\rangle.$$
The right hand side of the inequality above is a linear function. Hence the problem ($\ref{Merit}$) would be reduced to the following linear program:
\begin{eqnarray}\begin{split}\label{Merit1}
\mbox{Minimize } &\langle \nabla F_{\epsilon}(u),x\rangle&\\
\mbox {s.t.  } & Ax=b.&\\
\end{split}\end{eqnarray}
So in our iterative reweighted algorithm, we solve at step $l$ the following optimization problem:
\begin{eqnarray}\begin{split}\label{Merit2}
\mbox{Minimize } &\langle\nabla F_{\epsilon}(x^l),x\rangle&\\
\mbox {s.t.  } & Ax=b,&\\
\end{split}\end{eqnarray}
where $x^l$ is the solution of the previous iteration, and $\nabla F_{\epsilon}(x^l)$ are the weights. The reweighted $l_1$-algorithm can be defined as follows:

\begin{itemize}
\item Set $l$ as an index which counts the iterations, and choose a small enough $\epsilon>0$.
\item Step 0: Choose a starting point $x^{1}$. This can be obtained by solving the $l_1$-minimization problem.
\item Step $l$: Set $\omega^{l}=\nabla F_{\epsilon}(x^l)$, and Solve
\begin{equation}  x^{l+1}=\mbox{argmin}\{\langle \omega^{l},x\rangle:Ax=b\}.\end{equation}
\item Step $l+1$: If some termination criteria holds, stop. Otherwise, set $l\leftarrow l+1$, and go to step $l$.
\end{itemize}

An additional step can be added to the above algorithm concerning the choice of $\epsilon$. In this paper our updating rule is $\epsilon_{l+1}=0.5\epsilon_l$. In CWB algorithm, $\epsilon$ is updated as $\epsilon_{l+1}=\max\{|x^l|_{(i_0)},0.001\}$, where $i_0=\frac{m}{[4\log(\frac{n}{m})]}$, and $|x|_{i_0}$ is the biggest $i_0$ elements of $x$.\\

It is quit challenging to prove that under a mild condition, the reweighted $l_1$-algorithm converges to the sparsest solution of problem $(\ref{CMP})$. This is still an open questions in this field. However some progress have been made in this area \cite{1,57,26,24}. Mangasarian \cite{57} introduced a successive linearization algorithm(SLA) to find the solution of general complementarity problems, and proved that SLA algorithm terminates in finite number of iterations, and creates decreasing objective function values at each iteration. Furthermore he proved these values converge to a stationary point. Chen and Zhou in \cite{26} proved that the sequence generated by reweighted $l_1$-algorithm converges to a stationary point of a kind of truncated $l_p$-minimization problem ($0<p<1$). Similar results can also be found in \cite{22}. Recently, Zhao and Li \cite{1} defined a range space property(RSP) for matrices, under which he proved that the reweighted $l_1$-algorithm converges to certain sparse solution of the problem.

Following the framework of the reweighted $l_1$-algorithm in \cite{1}, we discuss some new merit functions. Before we go ahead, let's consider the following merit function
\begin{equation}\label{merit1} F_{\epsilon}(x)=\sum_{i=1}^n \log(|x_i|+\epsilon)+\sum_{i=1}^n(|x_i|+\epsilon)^p,\end{equation}
where $0<p<1$, which is mentioned in \cite{1}, based on which we will construct new merit functions. To verify that the above function is a merit function, one should check all the defined properties are satisfied. First, let's verify that this function is an approximation of $l_0$-norm function.\\

Indeed, it is easy to check that
$$\lim_{\epsilon\rightarrow 0} \left (n-\frac{\sum_{i=1}^n\log(|x_i|+\epsilon)+\sum_{i=1}^n(|x_i|+\epsilon)^p}{\log\epsilon}\right )=\|x\|_0.$$
Note that
$$\lim_{x\rightarrow \infty}F_{\epsilon}(x)=\infty,$$
which means that the function is coercive. It is clear that $F_{\epsilon}(x)=F_{\epsilon}(|x|)$, and the function is increasing.
In $R_+^n$, we have
$$\nabla F_{\epsilon}(x)=\left(\frac{1+(x_1+\epsilon)^p p}{x_1+\epsilon},...,\frac{1+(x_n+\epsilon)^p p}{x_n+\epsilon}\right)^T.$$
Also, for every $i=1,...,n$, we have
$$\lim_{(x_i,\epsilon)\rightarrow (0,0)}[ \nabla F_{\epsilon}(|x|)]_i=\lim_{(x_i,\epsilon)\rightarrow (0,0)} \frac{1+(|x_i|+\epsilon)^p p}{|x_i|+\epsilon}=\infty,~i=1,...,n,$$
and $\nabla F_{\epsilon}(x)$ is bounded when $\epsilon\rightarrow 0$. Since $0<p<1$ and $x_i>0$, we have $$ p(x_i+\epsilon)^p>p^2(x_i+\epsilon)^p,~i=1,...,n,$$ so
$$\frac{-1+(x_i+\epsilon)^pp^2-(x_i+\epsilon)^pp}{x_i+\epsilon}<0,~i=1,...,n,$$
and hence
$$\nabla^2F_{\epsilon}(x)=\mbox{diag}\left(\frac{-1+(x_i+\epsilon)^pp^2-(x_i+\epsilon)^pp}{x_i+\epsilon}\right)\prec 0,~i=1,...,n.$$
As seen, in $R_+^n$ the Hessian of the above merit function is negative definite, so the function  $F_{\epsilon}(x)$ is strictly concave.
From the above discussion one can define the following weights for the reweighted $l_1$-algorithm:
$$\omega_i=[ \nabla F_{\epsilon}(|x|)]_i=\frac{1+(|x_i|+\epsilon)^p p}{|x_i|+\epsilon},~i=1,...,n.$$
Note that the item (2) of the definition of a merit function implies that $[\nabla F_{\epsilon}(x^l)]_i\rightarrow \infty$ as $(x_i^l,\epsilon)\rightarrow (0,0)$, which means larger penalties(weights) for the smaller elements of $x$, at each iteration.\\

Now, we start to define a new merit function as follows
\begin{equation}\label{merit2}F_{\epsilon}(x)=\sum_{i=1}^n \log(\log(|x_i|+\epsilon+(|x_i|+\epsilon)^p)).\end{equation}
We verify this function is a merit function. Clearly, this function is an approximation of $l_0$-norm function. Because
$$\lim_{x\rightarrow 0}\frac{\log(\log(x_i+\epsilon+(x_i+\epsilon)^p))}{\log(\log(\epsilon))}=0,~\mbox{for}~ x_i\neq 0,$$
and
$$\lim_{x\rightarrow 0}\frac{\log(\log(x_i+\epsilon+(x_i+\epsilon)^p))}{\log(\log(\epsilon))}=1,~\mbox{for}~ x_i=0,$$
we conclude that
$$\lim_{x\rightarrow 0}\frac{\sum_{i=1}^n \log(\log(|x_i|+\epsilon+(|x_i|+\epsilon)^p))}{\log(\log(\epsilon))}=n-\|x\|_0.$$
In $R_+^n$, the gradient of $F_{\epsilon}(x)$ is given by
$$[\nabla F_{\epsilon}(x)]_i=\frac{1+\frac{(x_i+\epsilon)^pp}{x_i+\epsilon}}{(x_i+\epsilon+(x_i+\epsilon)^p)(\log(x_i+\epsilon+(x_i+\epsilon)^p))},$$
and since $\lim_{x_i\rightarrow0} x_i\log(x_i)=0$, we have
$$ \lim_{(x_i,\epsilon)\rightarrow (0,0)} [\nabla F_{\epsilon}(x)]_i=\infty.$$
Note that for every $x_i>0$,
$$\lim_{\epsilon\rightarrow 0}[\nabla F_{\epsilon}(x)]_i=\frac{1+px_i^{p-1}}{(x_i+x_i^p)\log(x_i+x_i^p)}=c_i,~i=1,...,n,$$
where $c_i$ is positive and bounded for every $i=1,...,n$. Also in $R_+^n$, $\nabla^2 F_{\epsilon}(x)$ is a diagonal matrix with the following entries on its diagonal,
$$[\nabla^2F_{\epsilon}(x)]_{ii}=\frac{\frac{(x_i+\epsilon)^pp^2}{(x_i+\epsilon)^2}-\frac{(x_i+\epsilon)^pp}{(x_i+\epsilon)^2}}{(x_i+\epsilon+(x_i+\epsilon)^p)\log(x_i+\epsilon+(x_i+\epsilon)^p)}$$
$$-\frac{\left(1+\frac{(x_i+\epsilon)^pp}{x_i+\epsilon}\right)^2}{\left(x_i+\epsilon+(x_i+\epsilon)^p\right)^2\log(x_i+\epsilon+(x_i+\epsilon)^p)}$$
$$-\frac{\left(1+\frac{(x_i+\epsilon)^pp}{x_i+\epsilon}\right)^2}{\left(x_i+\epsilon+(x_i+\epsilon)^p\right)^2\log\left(x_i+\epsilon+(x_i+\epsilon)^p\right)^2},~i=1,...,n.$$
Since, for every $i=1,...,n$, $[\nabla^2F_{\epsilon}(x)]_{ii}<0$, we have
$$ \nabla^2 F_{\epsilon}(x)\prec 0.$$
So the function $(\ref{merit1})$ is strictly concave, and it is a merit function.

Based on $(\ref{merit1})$, the reweighted $l_1$-algorithm choose the following weights:
$$\omega_i=\frac{1+\frac{(|x_i|+\epsilon)^pp}{|x_i|+\epsilon}}{(|x_i|+\epsilon+(|x_i|+\epsilon)^p)(\log(|x_i|+\epsilon+(|x_i|+\epsilon)^p))},~i=1,...,n.$$

In this paper, we refer to $W_1$ as the reweighted algorithm with the above weights. The Figure (15) shows the probability of success of $W_1$ algorithm via different choices of $\epsilon$. This figure demonstrates that $\epsilon=0.01$ works very good to locate the exact solution of the problem $(\ref{CMP})$, when the sparsity is 15. Clearly the above weights are related to the parameter $p$, so we tested the performance of $W_1$ algorithm for different sparsity of the solution, i.e, $k=5,10,15,20$, via different choices of $p$. Thirteen different values of $p$ have been tested (matrix $A$ has been normally distributed), and the result is summarized in Figure (14). Obviously the probability of success is higher when the sparsity of the solution is lower. This can be seen in Figure (14).\\

Another new merit function can be defined as follows
\begin{equation}\label{Merit2} F_{\epsilon}(x)=\frac{1}{p}\sum_{i=1}^n \left(\log(|x_i|+\epsilon+(|x_i|+\epsilon)^q)\right)^p,\end{equation}
where $0<p,q<1$, which is an approximation of $\|x\|_0$. In fact
$$n-\lim_{\epsilon\rightarrow0}\frac{\frac{1}{p}\sum_{i=1}^n \left(\log(|x_i|+\epsilon+(|x_i|+\epsilon)^q)\right)^p}{\frac{1}{p}\left(\log(\epsilon+\epsilon^q)\right)^p}=\|x\|_0.$$
In $R_+^n$, the gradient of $F_{\epsilon}(x)$ is given by
$$[\nabla F_{\epsilon}(x)]_i=\frac{\log\left(x_i+\epsilon+(x_i+\epsilon)^q\right)^p\left(1+\frac{(x_i+\epsilon)^qq}{x_i+\epsilon}\right)}{\left(x_i+\epsilon+(x_i+\epsilon)^q\right)\log\left(x_i+\epsilon+(x_i+\epsilon)^q\right)}.$$
Note that $\lim_{(x_i,\epsilon)\rightarrow (0,0)}[\nabla F_{\epsilon}(x)]_i=\infty$, and $\lim_{\epsilon\rightarrow 0}[\nabla F_{\epsilon}(x)]_i$ is bounded, for every fixed $x>0$.
In $R_+^n$, the Hessian is a diagonal matrix with the following diagonal elements
$$[\nabla^2F_{\epsilon}(x)]_{ii}=\frac{\log\left(x_i+\epsilon+(x_i+\epsilon)^p\right)^q}{\left(x_i+\epsilon+(x_i+\epsilon)^q\right)\log\left(x_i+\epsilon+(x_i+\epsilon)^q\right)}$$$$.\left(\frac{(p-1)\left(1+\frac{(x_i+\epsilon)^qq}{x_i+\epsilon}\right)^2}{\left(x_i+\epsilon+(x_i+\epsilon)^q\right)\log\left(x_i+\epsilon+(x_i+\epsilon)^q\right)}\right.+$$
$$\left.\frac{(x_i+\epsilon)^qq^2-(x_i+\epsilon)^qq}{(x_i+\epsilon)^2}-\frac{\left(1+\frac{(x_i+\epsilon)^qq}{x_i+\epsilon}\right)^2}{x_i+\epsilon+(x_i+\epsilon)^q}\right),~i=1,...,n.$$
where $0<p,q<1$. Clearly, $\nabla^2F_{\epsilon}(x)\prec0$, which implies the function is strictly concave. Thus, the function ($\ref{Merit2}$) is a merit function, so we may choose the following weights in our algorithm:
$$\omega_i=\frac{\log\left(|x_i|+\epsilon+(|x_i|+\epsilon)^q\right)^p\left(1+\frac{(|x_i|+\epsilon)^qq}{|x_i|+\epsilon}\right)}{\left(|x_i|+\epsilon+(|x_i|+\epsilon)^q\right)\log\left(|x_i|+\epsilon+(|x_i|+\epsilon)^q\right)},~i=1,...,n.$$

We refer to $W_2$ as the reweighted algorithm with the weights above. The Figures (11),(12),(13) show the performance of $W_2$ algorithm for finding the exact solution of the problem $(\ref{CMP})$  for different choices of the weights parameters, $p$ and $q$, and for different fixed sparsity of the solution, i.e., $k=5,10,15,20$.\\

\textbf{Remark.}We see from above that the $\log$ function plays a vital rule in constructing a merit function. As pointed in \cite{1,62}, the $\log$ function can enhance the concavity of a given function without affecting its coercivity and monotonicity. For the convergency analysis of the reweighted $l_1$-algorithms based on the class of merit functions that defined at the beginning of this chapter, one may refer to the Theorems 3.9 and 3.11 in \cite{1}, where it has been shown that under the so-called RSP condition, the algorithm may converge to a solution of problem ($\ref{CMP}$) with certain level of sparsity.\\

\section{Numerical Experiments}
In this section, we compare the performance of the algorithms above for finding the exact solution of the problem $(\ref{CMP})$ through the numerical tests. We compare the following algorithms in our numerical experiments.\\
$l_1$-min:
\begin{eqnarray}\begin{split}\label{l1}
\mbox{Minimize } &\|x\|_1& \\
\mbox {s.t.  } & Ax=b,&\\
\end{split}\end{eqnarray}
CWB(Candes, Wakin, Boyd):
\begin{eqnarray}\begin{split}\label{CWB}
x^{l+1}=\mbox{argmin } &\sum_{i=1}^n\frac{1}{|x_i^l|+\epsilon^l}|x_i|& \\
\mbox {s.t.  } & Ax=b,&\\
\end{split}\end{eqnarray}
$W_1$:
\begin{eqnarray}\begin{split}\label{LW1}
x^{l+1}=\mbox{argmin } &\sum_{i=1}^n\frac{1+\frac{(|x_i^l|+\epsilon^l)^pp}{|x_i^l|+\epsilon^l}}{(|x_i^l|+\epsilon^l+(|x_i^l|+\epsilon^l)^p)(\log(|x_i^l|+\epsilon^l+(|x_i^l|+\epsilon^l)^p))}|x_i|& \\
\mbox {s.t.  } & Ax=b,&\\
\end{split}\end{eqnarray}
$W_2$:
\begin{eqnarray}\begin{split}\label{LW2}
x^{l+1}=\mbox{argmin } &\sum_{i=1}^n\frac{\log\left(|x_i^l|+\epsilon^l+(|x_i^l|+\epsilon^l)^q\right)^p\left(1+\frac{(|x_i^l|+\epsilon^l)^qq}{|x_i^l|+\epsilon^l}\right)}{\left(|x_i^l|+\epsilon+(|x_i^l|+\epsilon^l)^q\right)\log\left(|x_i^l|+\epsilon^l+(|x_i^l|+\epsilon^l)^q\right)}|x_i|& \\
\mbox {s.t.  } & Ax=b,&\\
\end{split}\end{eqnarray}
where $0<p,q<1$, $A\in \mathbb{R}^{50\times 200},~b\in \mathbb{R}^{50},~\mbox{and}~x\in\mathbb{R}^{200}$.\\

In our numerical works, we randomly generated the matrix $A\in \mathbb{R}^{50\times 200}$, and for a fixed sparsity, we randomly generated the solution vector $x\in \mathbb{R}^{200}$. We tested 100 randomly generated marices, $A$, for different level of $k$-sparsity of the solution, i.e., $k=1,2,...,26$. The matrix $A$ (the problem data) was randomly generated based on different statistical distributions. Most of the previous numerical experiments in the literature usually use normally distributed matrices.\\

The distributions that we considered were Normal $(N(\mu, \sigma))$ with the parameters $\mu=0$ and $\sigma=1$, Poisson $(Pois(\lambda))$ with the parameter $\lambda=2$, Exponential $(Exp(\mu))$ with the parameter $\mu=5$, F-distribution $(F(\alpha,\beta))$ with the parameters $\alpha=1$ and $\beta=6$, Gamma distribution $(Gam(a,b))$ with parameters $a=5$ and $b=10$, and Uniform distribution $(U(N))$ with the parameter $N=10$. The probability of success of the 4 algorithms mentioned above, i.e, $l_1$-min, $CWB$, $W_1$, $W_2$ have been compared via different sparsity of the solution, and through all the above differently distributed matrices $A$. On a laptop with a Core 2 Duo CPU (2.00 GHz, 2.00GHz) and 4.00 GB of RAM memory, each comparing figure took approximately 14-hours time (in average).\\

The updating rule $\epsilon^{l+1}=\frac{1}{2}\epsilon^l$ was used, at each iteration $l$. The choice of $\epsilon$ is crucial for reweighted $l_1$-algorithms. Hence, we have also tested the algorithms by applying Candes, Wakin, Boyd(CWB) updating rule for $\epsilon$, and also a fixed $\epsilon=0.01$. These figures demonstrate how these choices of $\epsilon$ may affect the performance of the algorithms.\\

As seen, the weights in $W_1$ and $W_2$ vary for different values of $p$ and $p,q$. Therefore, we have tried different choices of $p$ and $q$ to find out how they may affect the success probability for $W_1$ and $W_2$ algorithms.\\

In the Figure (3), the matrix A has been generated from $Exp(\mu)$, with $\mu=5$. We set $p=0.05$ in $W_1$, and $p=q=0.05$ in $W_2$. As shown, all of the algorithms are very successful when $\|x\|_0<7$. When $7<\|x\|_0<11$, CWB, $W_1$ and $W_2$ almost perform the same as each other, but when $\|x\|_0>11$, $W_1$ and $W_2$ outperform the CWB algorithm. All of the algorithms fail when the cardinality of the solution is above 25, i.e, $\|x\|_0>25$.\\

In Figure (4), the matrix A has been generated from $Exp(\mu)$, with $\mu=5$, as in Figure (3). However, in this case we used different values for $p$ and $q$. We chose $p=q=0.4$, which is much larger than $0.05$. As expected, both $W_1$ and $W_2$ perform significantly worse than the case of $p=q=0.05$. Even for lower sparsity, both algorithms fail to locate the exact sparse solution with a high probability.\\

In Figure (5), the matrix $A$ has been generated from $F(\alpha, \beta)$, with $\alpha=1$ and $\beta=6$, and we set $p=q=0.05$.  As shown, all of the algorithms start failing when the cardinality of the solution is higher than 4, i.e, $\|x\|_0>4$. $W_1$ and $W_2$ perform better than $CWB$ for higher cardinality of the solution, and $W_2$ is slightly better than $W_1$ in general.\\

In Figure (6), the matrix $A$ has been generated from $Gam(a,b)$, with $a=5$ and $b=10$, and we set $p=q=0.05$. For lower cardinality of the solution, $CWB$ and $W_1$ perform slightly better than $W_2$ when $\|x\|_0<8$. Also $CWB$, $W_1$, and $l_1$-min are completely successful for locating the exact solution, when $\|x\|_0<8$. But for $8<\|x\|_0<10$, only $CWB$ and $W_1$ are successful. $l_1$-min, $CWB$, $W_1$ and $W_2$ fail when $\|x\|_0>17$, $\|x\|_0>21$, $\|x\|_0>24$, $\|x\|_0>26$, respectively. Therefore $W_1$ and $W_2$ perform significantly better for higher cardinality of the solution.\\

In Figure (7), the matrix $A$ has been generated from $N(\mu,\sigma)$, with $\mu=0$ and $\sigma=1$, and we set $p=q=0.05$. As shown, $l_1$-min, $CWB$, and $W_1$ are very successful for finding the sparsest solution of the system when $\|x\|_0<8$. $W_1$ and $W_2$ perform better than the other two algorithms for higher cardinality of the solution.\\

In Figure (8), the matrix $A$ has been generated from $N(\mu,\sigma)$, with $\mu=0$ and $\sigma=1$ as in the Figure(7). However, we chose bigger values for $p$ and $q$, i.e, $p=q=0.4$. For large values of $p$ and $q$, $W_2$ starts failing for $\|x\|_0>4$, and performs much worst than $W_1$, $CWB$, and $l_1$-min. Also, for higher cardinality of the solution $CWB$ performs better than $W_1$ and $W_2$. Hence, from this figure and the Figure (4), one may conclude that smaller values for $p$ and $q$ should be chosen in order to achieve better results. Note that for large values of $p$ and $q$ the merit functions in $W_1$ and $W_2$ are not good concave approximations of $l_0$-norm.\\

In Figure (9), the matrix $A$ has been generated form $U(N)$, with $N=10$, and we set $p=q=0.05$. All of the algorithms except $W_2$ are successful for finding the sparsest solution of the system when $\|x\|_0<9$. For $9<\|x\|_0<12$, $CWB$ performs slightly better than $W_1$ and $W_2$. But for higher cardinality of the solution, $W_1$ and $W_2$ outperform $l_1$-min and $CWB$.\\

In Figure(10), the matrix $A$ has been generated from $Pois(\lambda)$, with $\lambda=5$, and we set $p=q=0.05$. All of the algorithms except $W_2$ are successful for finding the sparsest solution of the system when $\|x\|_0<9$. For higher cardinality of the solution, $W_1$ and $W_2$ outperform the other algorithms.\\

Clearly, for small values of $p$, the best algorithm is $W_1$ in general, i.e. for different cardinality of the solution and for different tested distributions. For all of the different tested distributions, both $W_1$ and $W_2$ (for small choices of $p$ and $q$) outperform $CWB$ when the cardinality of the solution is higher.\\

In the Figures (11), (12), (13), we focused on the performance of $W_2$ algorithms for different values of $p$ and $q$ via different fixed cardinality of the solution. In the Figure (11), we fixed $p=0.08$ and set different values of $q$. We examined the probability of success of $W_2$ for different fixed sparsity of 5,10,15,20. As expected, when cardinality of the solution is lower the success probability of $W_2$ is higher. As seen, the probability of success for fixed sparsity of 5 is the highest, and the probability of success for fixed sparsity of 20 is the lowest. The Figures (12) and (13) show the same results for fixed $p=0.4$ and $p=0.8$, respectively.\\

In the Figure (14), the performance of $W_1$ has been tested using different choices of $p$ and different fixed sparsity of the solution. As seen, in Figure (14), when $p$ increases from 0.04 to 1, the probability of success of the algorithm becomes lower(except some jumps). As shown, for different fixed sparsity of 5,10,15,20 the highest probability of success was achieved when $p=0.04$. Looking back to the merit function defined for the $W_1$ algorithm, one may see that for smaller values of $p$ the function is a better concave approximation of $l_0$-norm.\\

As we have discussed before, the choice of $\epsilon$ for the reweighted $l_1$-algorithm is important. Either very small or very big $\epsilon$ may result in improper weights, which may cause the failure of the algorithms. In Figure (15), we fixed the sparsity of the solution ($k=15$) and set $p=0.05$. Different choices of $\epsilon$ have been tested to suggest what $\epsilon$ might be the good one for which $W_1$ performs better. The matrix $A$ has been generated from $N(0,1)$. As shown, when $\epsilon$ tends to zero (e.g. $\epsilon\approx 0.0001$), or when $\epsilon$ is big (e.g. $\epsilon \approx 0.1$), the probability of success decreases. Our numerical experiments, in Figure (15), show that $\epsilon=0.01$ is a good choice for the weights in $W_1$ algorithm.\\

In Figure (16), we fixed $\epsilon=0.01$(with no updating rule) and compared the performance of $l_1$-min, $CWB$, $W_1$, $W_2$. Like Figure (7), the matrix has been generated from $N(0,1)$, and we set $p=q=0.05$. Our numerical experiment show that $W_1$ and $W_2$ significantly outperform $CWB$ especially for higher cardinality of the solution. Comparing Figure (16) and Figure (7), one may conclude that even for a fixed $\epsilon$, if chosen correctly, both $W_1$ and $W_2$ algorithms may perform very well to find a sparse solution.\\

Again to show that how important the choice of $\epsilon$ is, we compare the performance of $l_1$-min, $CWB$, $W_1$, $W_2$ based on the Candes updating rule. As seen in Figure(17), $CWB$ outperforms both $W_1$ and $W_2$ for lower cardinality of the solution, i.e, when $\|x\|_0<12$ in our numerical experiments.\\

\section{Conclusion}
We introduced a few concave approximations for the function $\|x\|_0$. These approximations can be employed to define new weights for the reweighted $l_1$-algorithms, which are used to locate the sparse solution of a linear system of equations. Through numerical experiments, we compared the performance of these reweighted algorithms and some existing reweighted algorithms when applied to linear systems with different statistically distributed matrices $A$, and with different sparsity of the solution. We have also explained when the new reweighted algorithms outperform some existing algorithms in different situations. We have also illustrated that how different choices of $p$ and $q$ may affect the performance of the algorithms. In addition, we have shown that the choices of $\epsilon$ may remarkably affect the performance of these algorithms as well.
~\\
~\\
~\\
~\\
~\\
~\\
~\\
~\\
~\\
~\\
~\\
~\\
~\\
~\\
~\\
~\\
~\\
~\\
~\\
~\\
~\\
~\\
~\\
~\\
~\\
~\\
~\\
~\\
~\\
~\\
~\\
~\\
~\\

\begin{figure}[htp]
\centering
$\begin{array}{c}
\includegraphics[width=0.5\textwidth,totalheight=0.25\textheight]{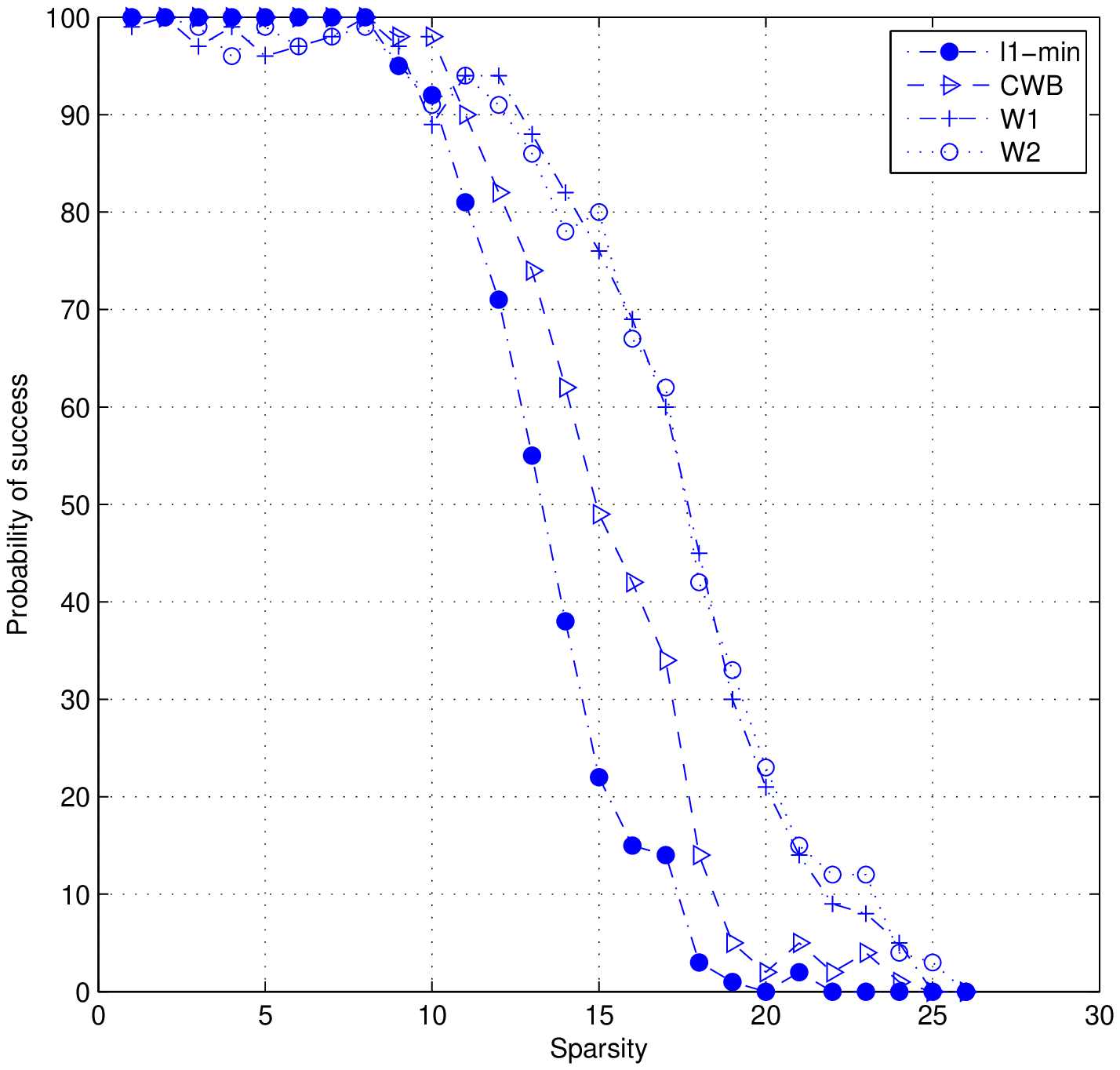}
\label{FIG2}
\end{array}$
\caption{Comparing the performance of $l_1$-min, CWB, $W_1$, $W_2$ minimization via the probability of success for finding the exact $k$-sparse solution of $Ax=b$, where $A\in \mathbb{R}^{50\times 200}$, $b\in \mathbb{R}^{50}$, $p=q=0.05$. Matrix $A$ has been generated from Exponential distribution. 100 randomly generated matrices have been tested for different sparsity of $k=1,...,26$.}
$\begin{array}{c}
\includegraphics[width=0.5\textwidth,totalheight=0.25\textheight]{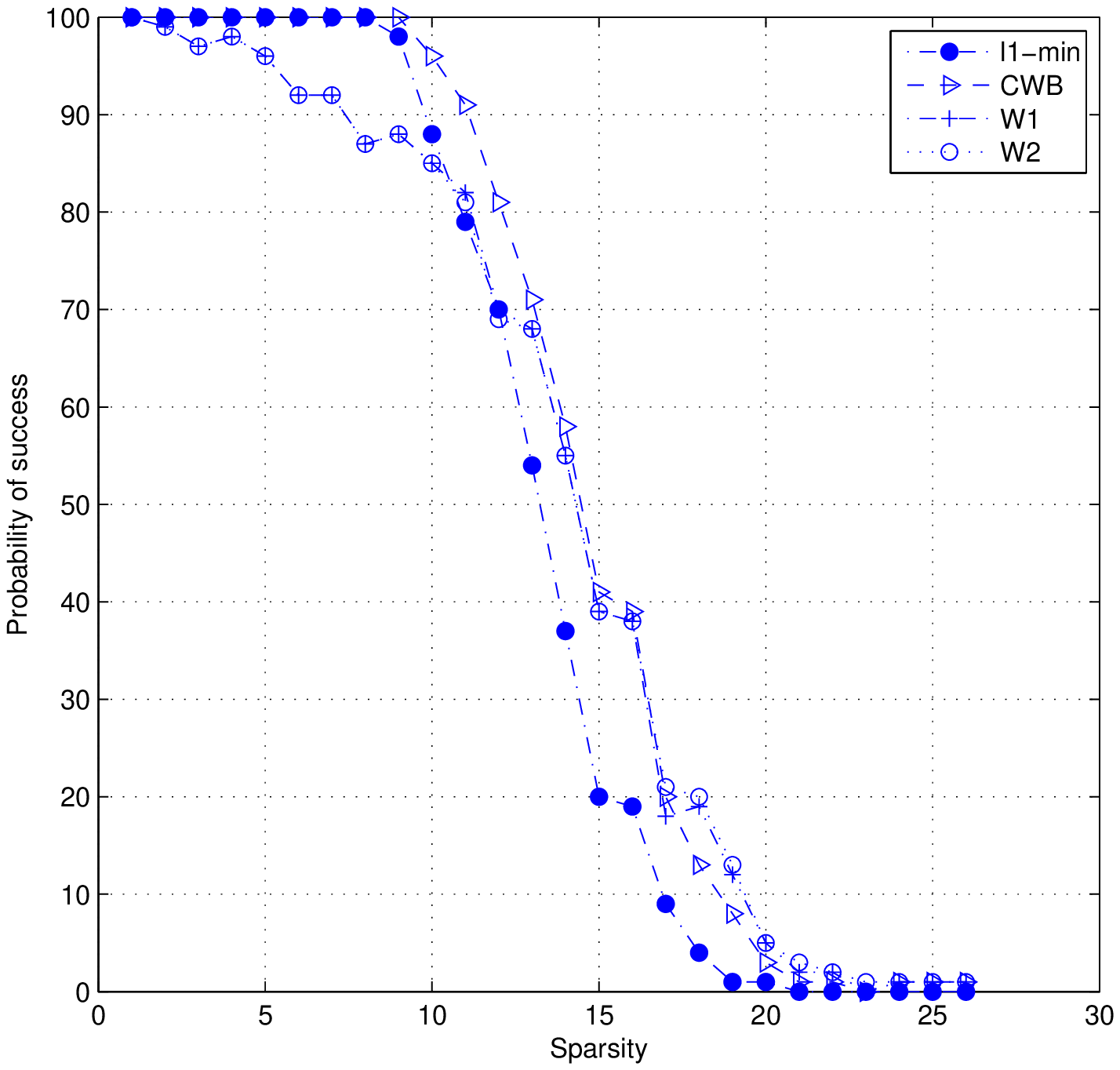}
\label{FIG3}
\end{array}$
\caption{Comparing the performance of $l_1$-min, CWB, $W_1$, $W_2$ minimization via the probability of success for finding the exact $k$-sparse solution of $Ax=b$, where $A\in \mathbb{R}^{50\times 200}$, $b\in \mathbb{R}^{50}$, $p=q=0.4$. Matrix $A$ has been generated from Exponential distribution. 100 randomly generated matrices have been tested for different sparsity of $k=1,...,26$.}
\end{figure}

\begin{figure}[htp]
\centering
$\begin{array}{c}
\includegraphics[width=0.5\textwidth,totalheight=0.25\textheight]{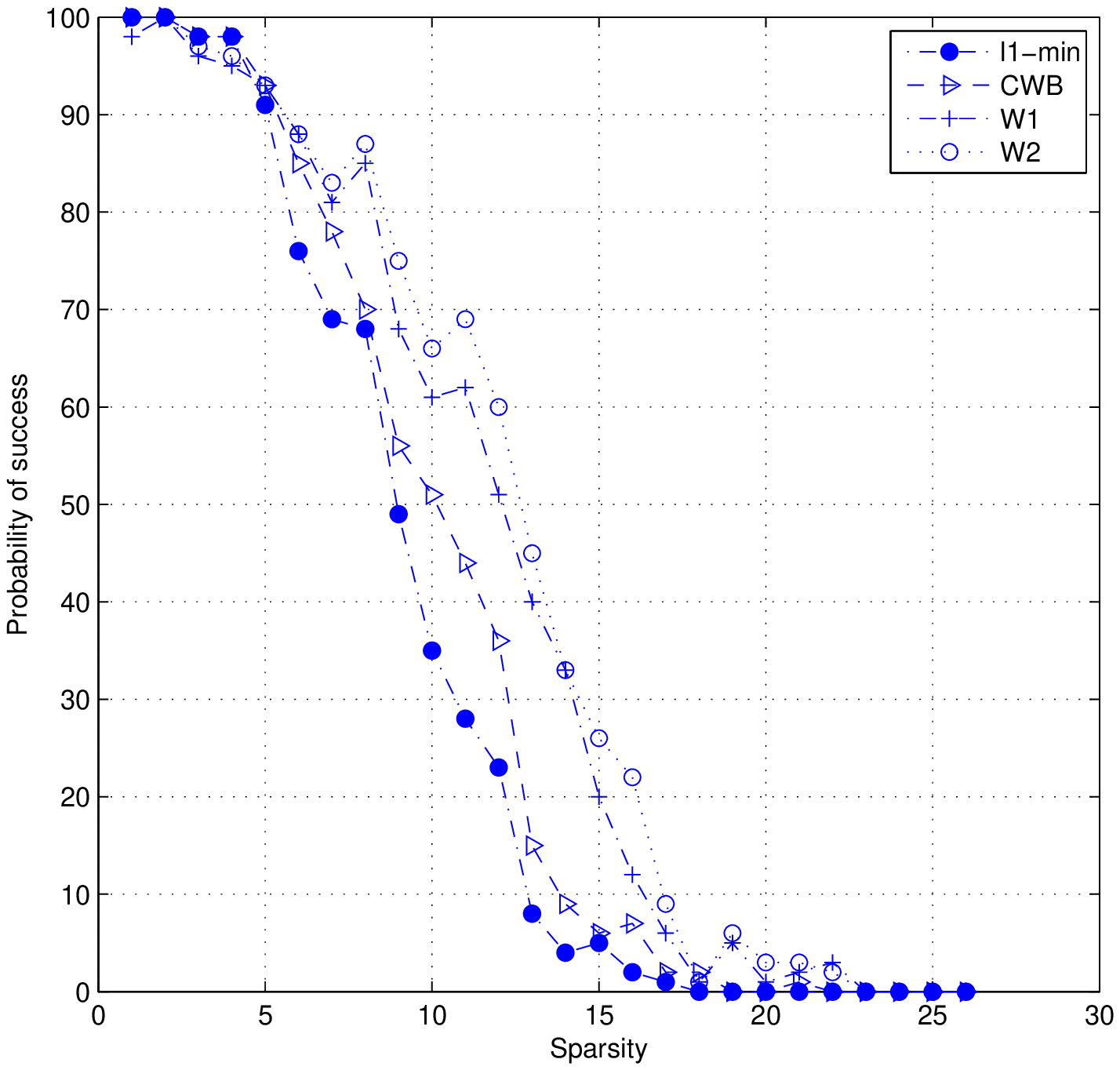}
\label{FIG4}
\end{array}$
\caption{Comparing the performance of $l_1$-min, CWB, $W_1$, $W_2$ minimization via the probability of success for finding the exact $k$-sparse solution of $Ax=b$, where $A\in \mathbb{R}^{50\times 200}$, $b\in \mathbb{R}^{50}$, $p=q=0.05$. Matrix $A$ has been generated from F-distribution. 100 randomly generated matrices have been tested for different sparsity of $k=1,...,26$.}
$\begin{array}{c}
\includegraphics[width=0.5\textwidth,totalheight=0.25\textheight]{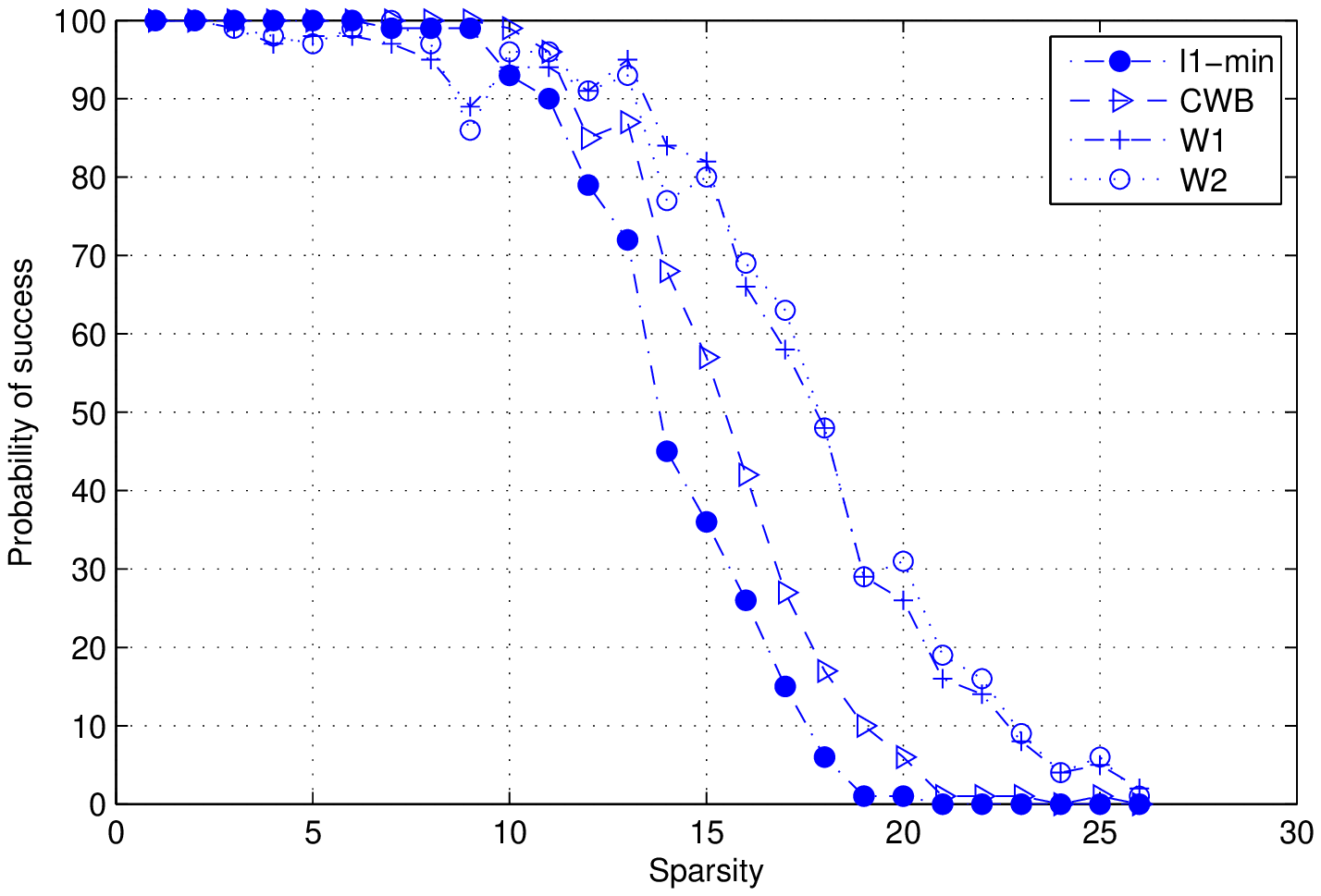}
\label{FIG5}
\end{array}$
\caption{Comparing the performance of $l_1$-min, CWB, $W_1$, $W_2$ minimization via the probability of success for finding the exact $k$-sparse solution of $Ax=b$, where $A\in \mathbb{R}^{50\times 200}$, $b\in \mathbb{R}^{50}$, $p=q=0.05$. Matrix $A$ has been generated from Gamma distribution. 100 randomly generated matrices have been tested for different sparsity of $k=1,...,26$.}
\end{figure}

\begin{figure}[htp]
\centering
$\begin{array}{c}
\includegraphics[width=0.5\textwidth,totalheight=0.25\textheight]{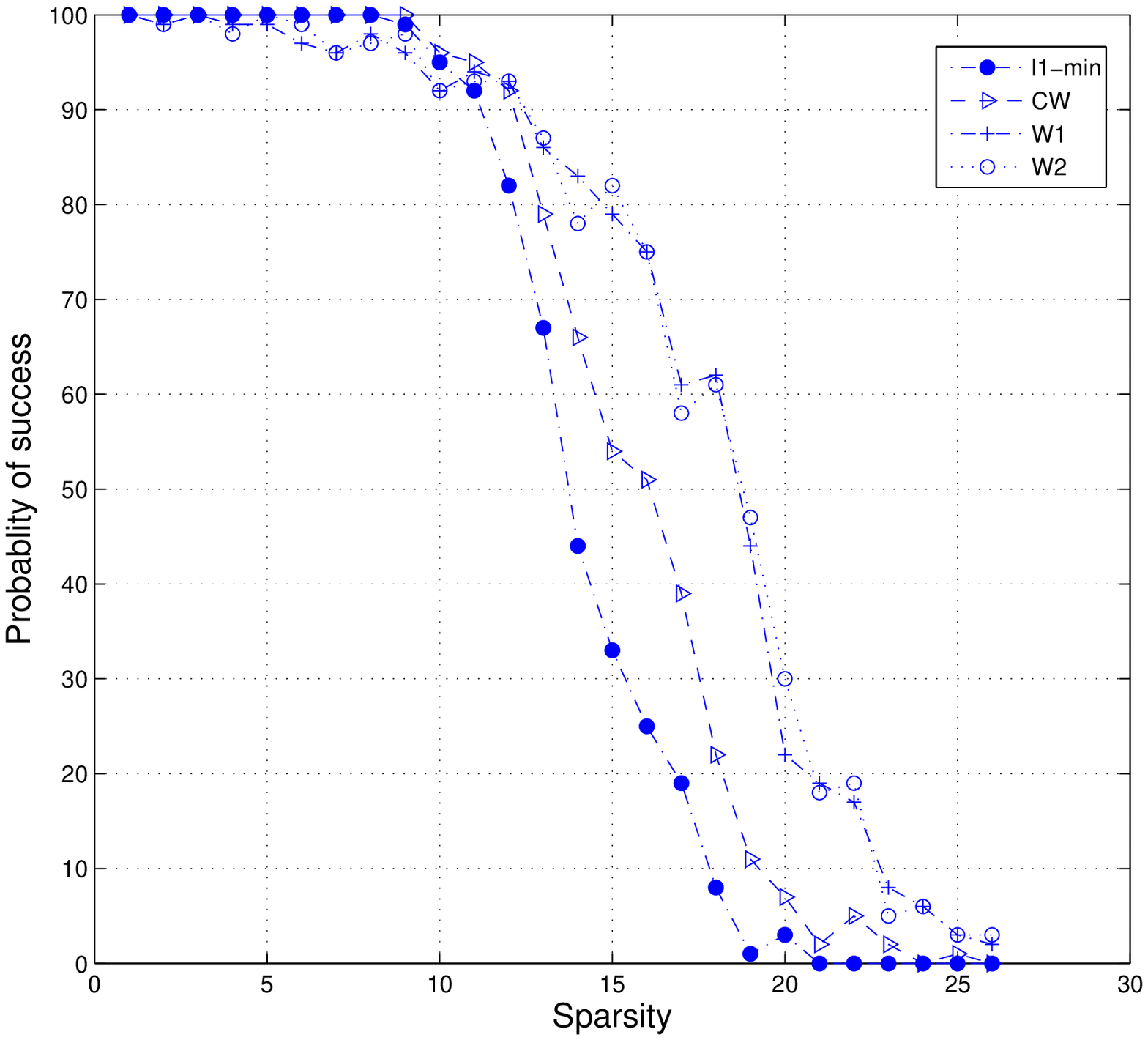}
\label{FIG6}
\end{array}$
\caption{Comparing the performance of $l_1$-min, CWB, $W_1$, $W_2$ minimization via the probability of success for finding the exact $k$-sparse solution of $Ax=b$, where $A\in \mathbb{R}^{50\times 200}$, $b\in \mathbb{R}^{50}$, $p=q=0.05$. Matrix $A$ has been generated from Normal distribution. 100 randomly generated matrices have been tested for different sparsity of $k=1,...,26$.}
$\begin{array}{c}
\includegraphics[width=0.5\textwidth,totalheight=0.25\textheight]{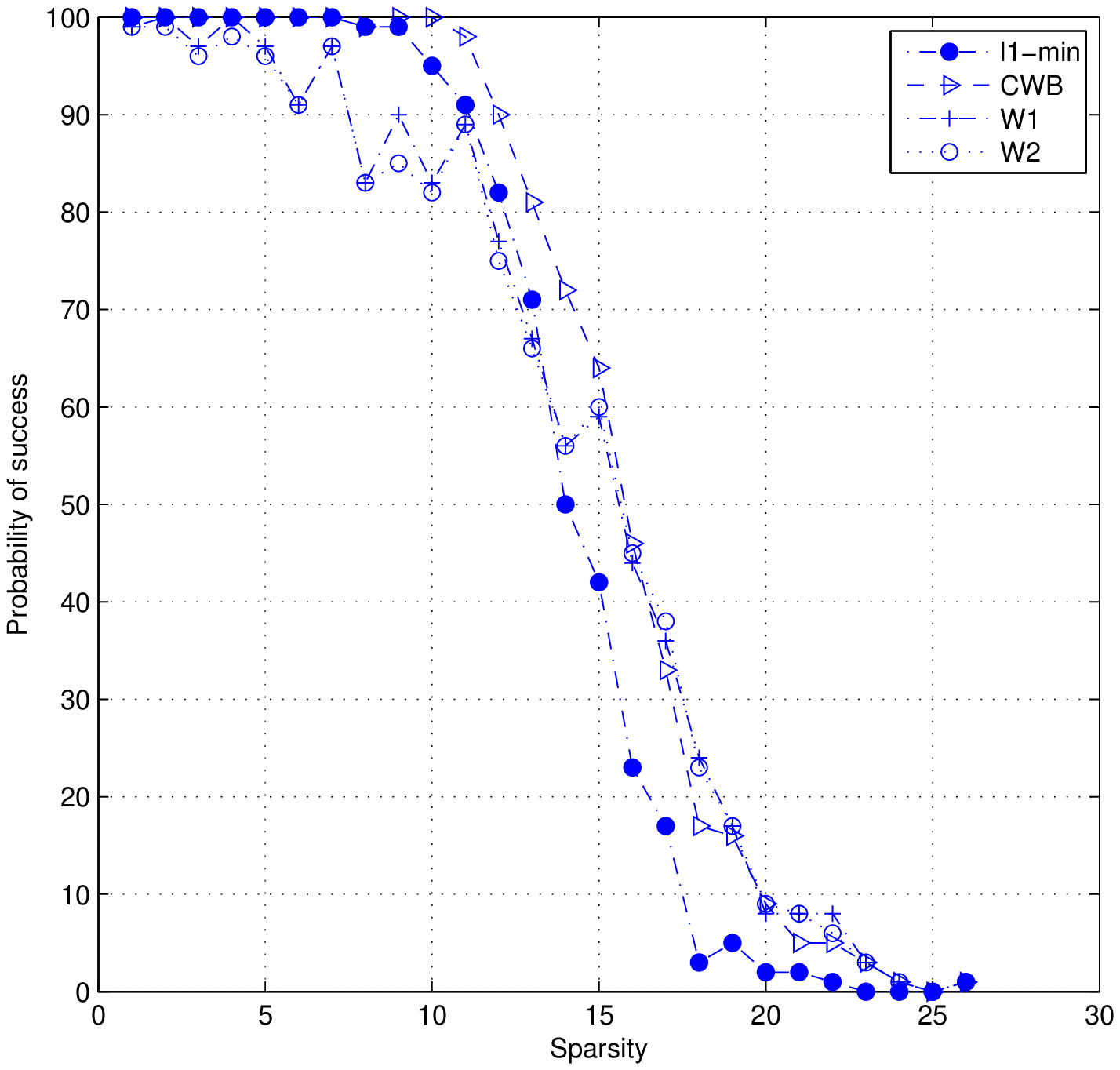}
\label{FIG7}
\end{array}$
\caption{Comparing the performance of $l_1$-min, CWB, $W_1$, $W_2$ minimization via the probability of success for finding the exact $k$-sparse solution of $Ax=b$, where $A\in \mathbb{R}^{50\times 200}$, $b\in \mathbb{R}^{50}$, $p=q=0.4$. Matrix $A$ has been generated from Normal distribution. 100 randomly generated matrices have been tested for different sparsity of $k=1,...,26$.}
\end{figure}

\begin{figure}[htp]
\centering
$\begin{array}{c}
\includegraphics[width=0.5\textwidth,totalheight=0.25\textheight]{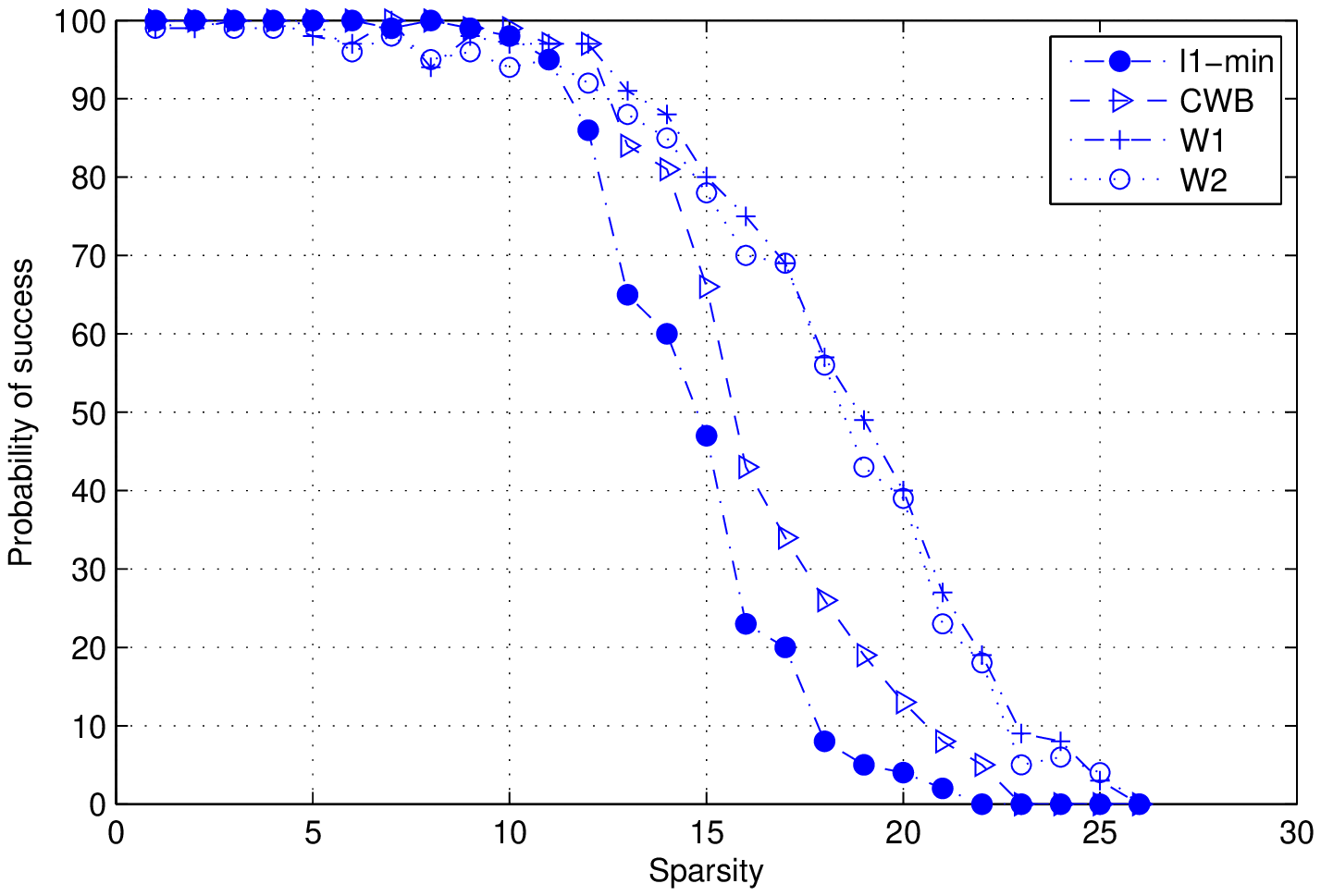}
\label{FIG8}
\end{array}$
\caption{Comparing the performance of $l_1$-min, CWB, $W_1$, $W_2$ minimization via the probability of success for finding the exact $k$-sparse solution of $Ax=b$, where $A\in \mathbb{R}^{50\times 200}$, $b\in \mathbb{R}^{50}$, $p=q=0.05$. Matrix $A$ has been generated from Uniform distribution. 100 randomly generated matrices have been tested for different sparsity of $k=1,...,26$.}
$\begin{array}{c}
\includegraphics[width=0.5\textwidth,totalheight=0.25\textheight]{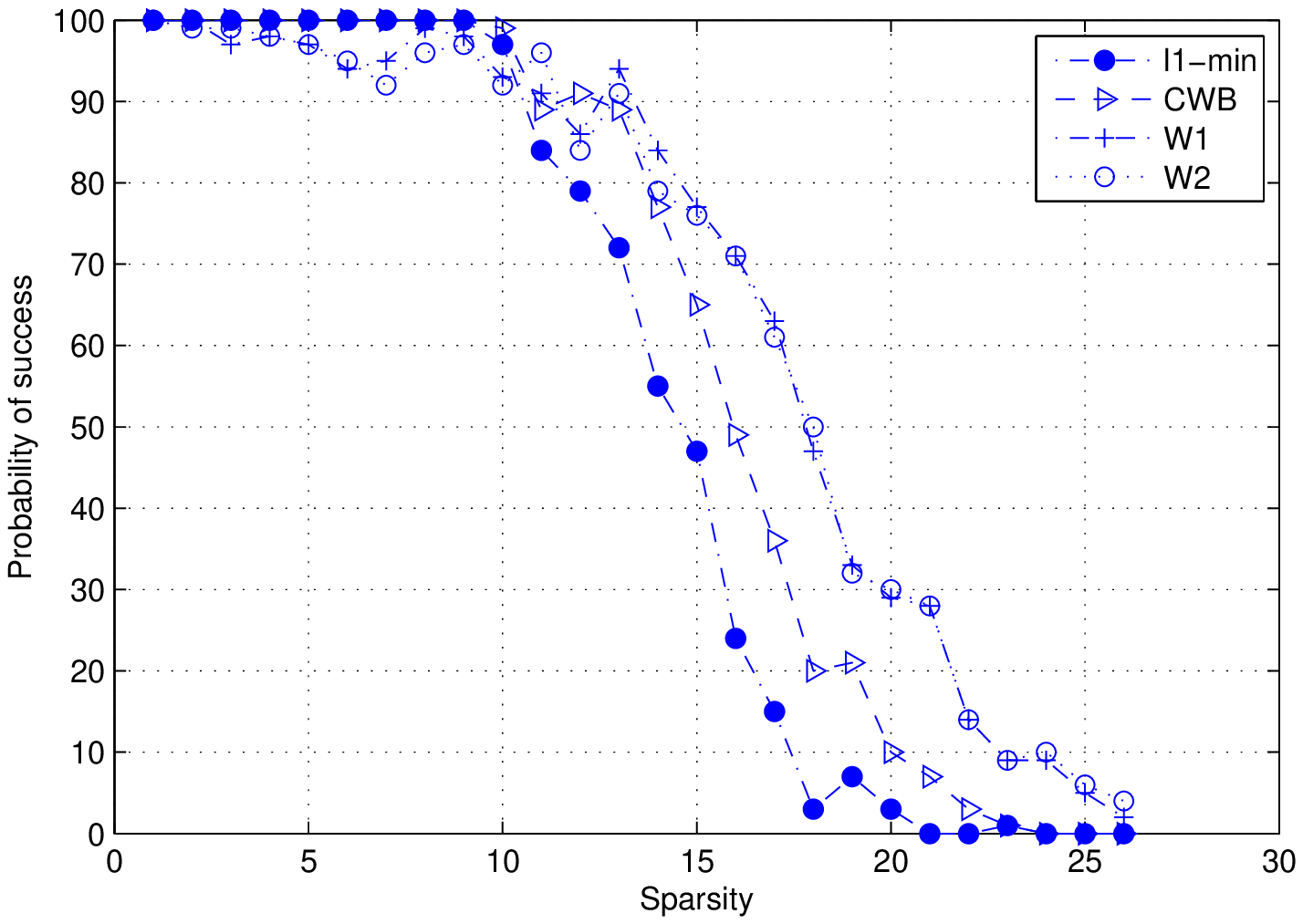}
\label{FIG9}
\end{array}$
\caption{Comparing the performance of $l_1$-min, CWB, $W_1$, $W_2$ minimization via the probability of success for finding the exact $k$-sparse solution of $Ax=b$, where $A\in \mathbb{R}^{50\times 200}$, $b\in \mathbb{R}^{50}$, $p=q=0.05$. Matrix $A$ has been generated from Poisson distribution. 100 randomly generated matrices have been tested for different sparsity of $k=1,...,26$.}
\end{figure}

\begin{figure}[htp]
\centering
$\begin{array}{c}
\includegraphics[width=0.5\textwidth,totalheight=0.25\textheight]{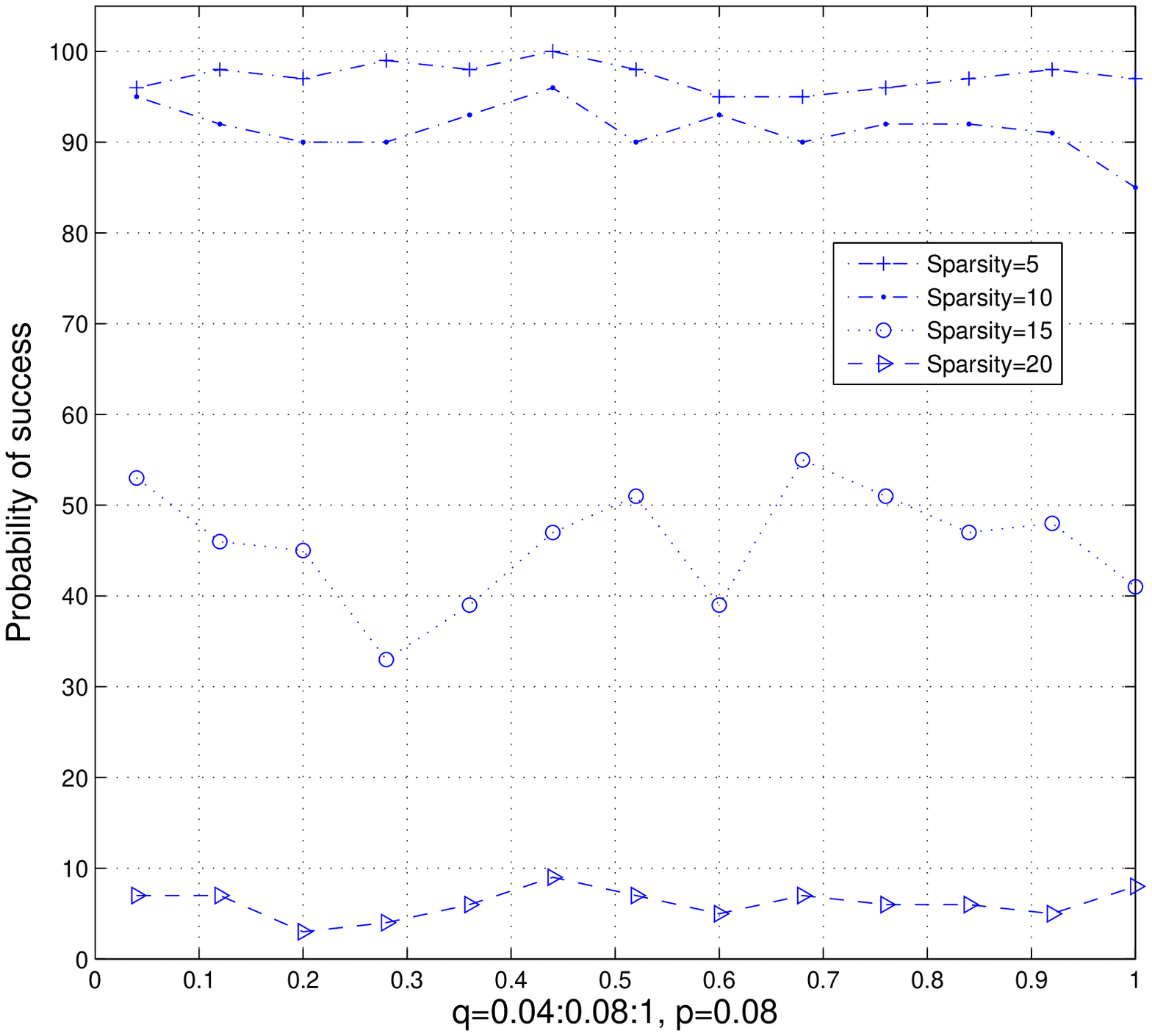}
\label{FIG10}
\end{array}$
\caption{Comparing the performance of $W_2$ minimization for different $q=0.04:0.08:1$, $p=0.08$ via the probability of success for finding the exact $k$-sparse solution of $Ax=b$, where $A\in \mathbb{R}^{50\times 200}$, $b\in \mathbb{R}^{50}$. Matrix $A$ has been generated from normal distribution. 100 randomly generated matrices have been tested for different sparsity of $k=5,10,15,20$.}
$\begin{array}{c}
\includegraphics[width=0.5\textwidth,totalheight=0.25\textheight]{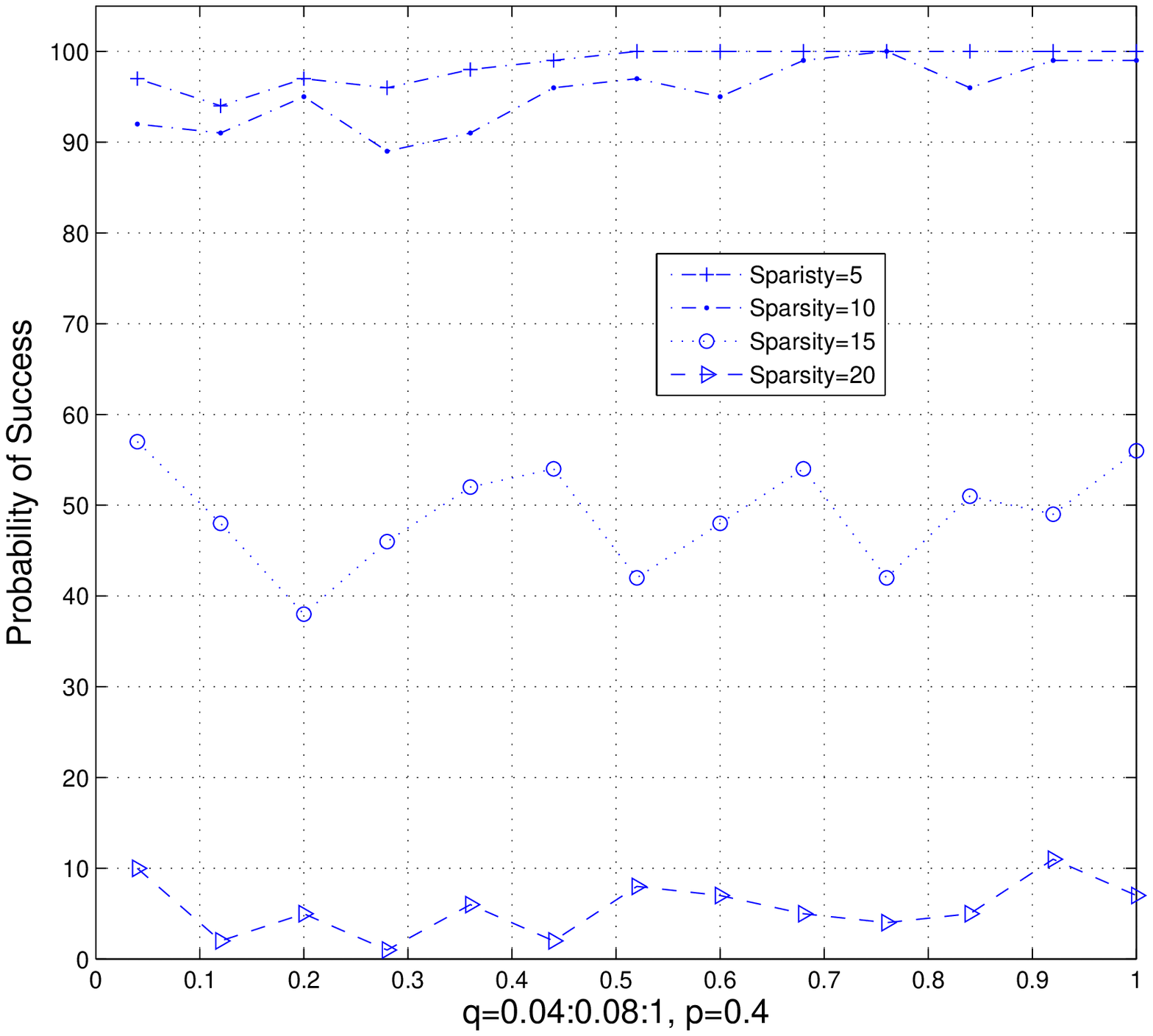}
\label{FIG11}
\end{array}$
\caption{Comparing the performance of $W_2$ minimization for different $q=0.04:0.08:1$, $p=0.4$ via the probability of success for finding the exact $k$-sparse solution of $Ax=b$, where $A\in \mathbb{R}^{50\times 200}$, $b\in \mathbb{R}^{50}$. Matrix $A$ has been generated from normal distribution. 100 randomly generated matrices have been tested for different sparsity of $k=5,10,15,20$.}
\end{figure}

\begin{figure}[htp]
\centering
$\begin{array}{c}
\includegraphics[width=0.5\textwidth,totalheight=0.25\textheight]{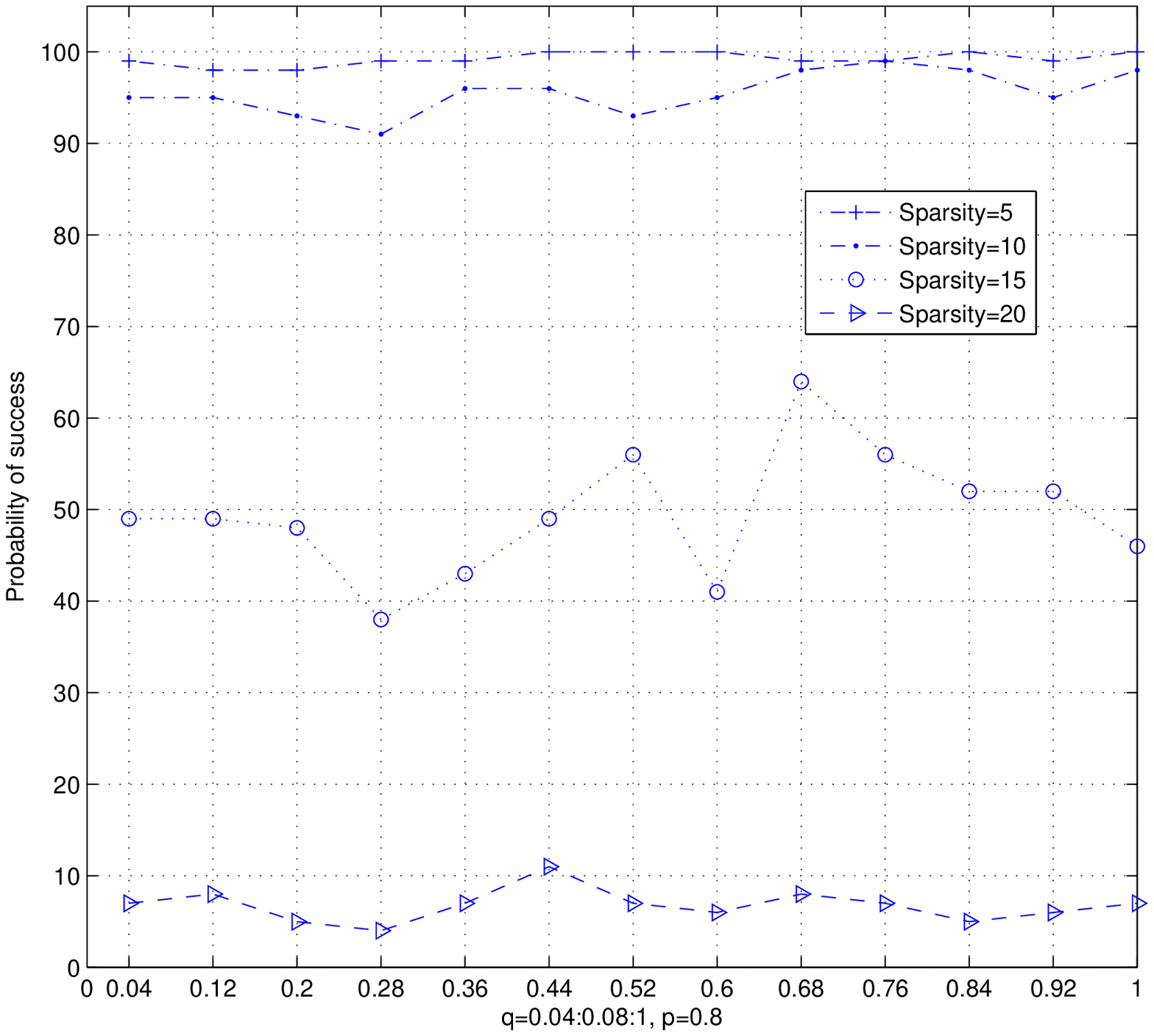}
\label{FIG12}
\end{array}$
\caption{Comparing the performance of $W_2$ minimization for different $q=0.04:0.08:1$, $p=0.8$ via the probability of success for finding the exact $k$-sparse solution of $Ax=b$, where $A\in \mathbb{R}^{50\times 200}$, $b\in \mathbb{R}^{50}$. Matrix $A$ has been generated from normal distribution. 100 randomly generated matrices have been tested for different sparsity of $k=5,10,15,20$.}
$\begin{array}{c}
\includegraphics[width=0.5\textwidth,totalheight=0.25\textheight]{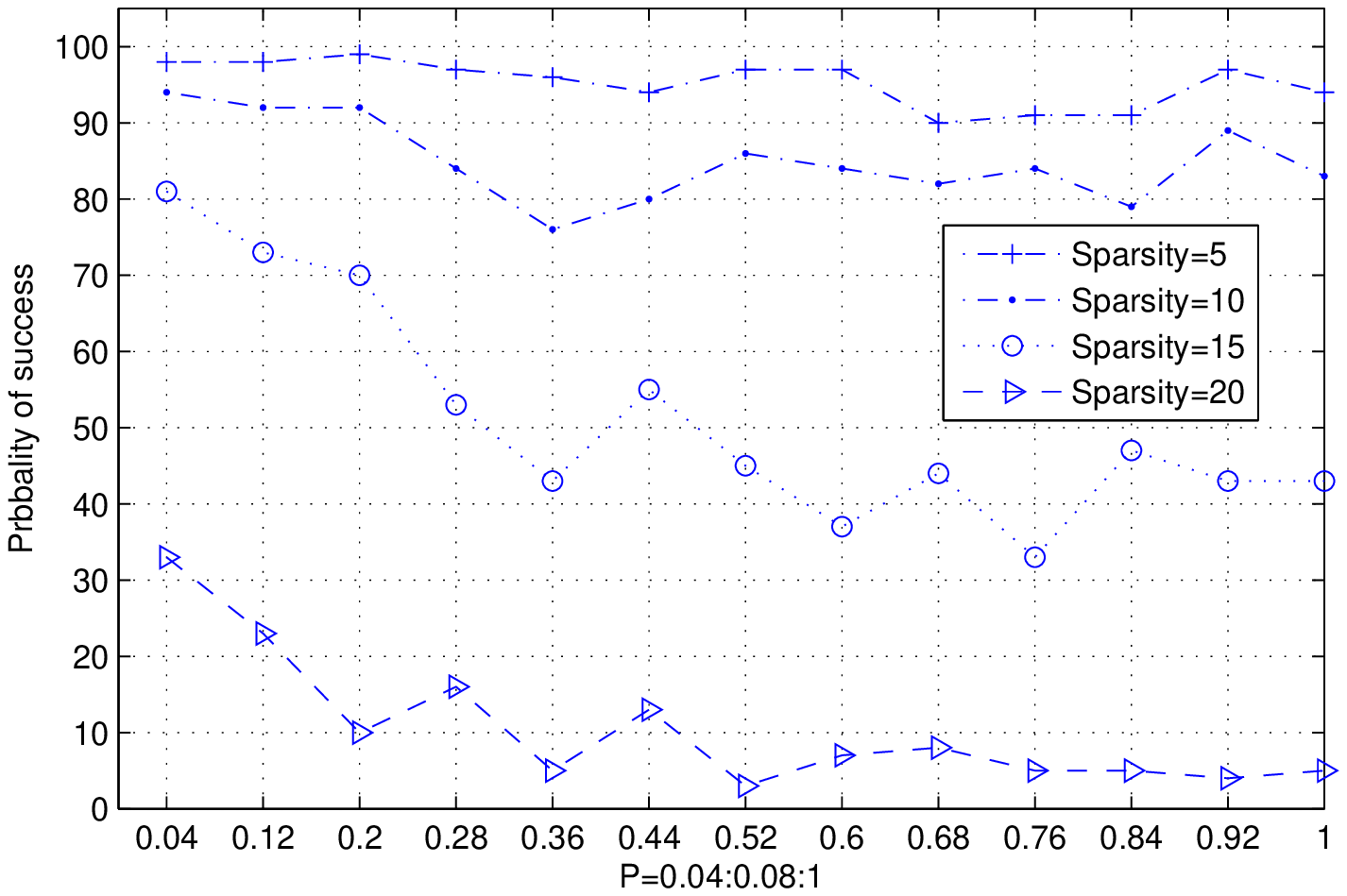}
\label{FIG9}
\end{array}$
\caption{Comparing the performance of $W_1$ minimization for different $p=0.04:0.08:1$ via the probability of success for finding the exact $k$-sparse solution of $Ax=b$, where $A\in \mathbb{R}^{50\times 200}$, $b\in \mathbb{R}^{50}$. Matrix $A$ has been generated from normal distribution. 100 randomly generated matrices have been tested for different sparsity of $k=5,10,15,20$.}
\end{figure}

\begin{figure}[htp]
\centering
$\begin{array}{c}
\includegraphics[width=0.5\textwidth, totalheight=.25\textheight]{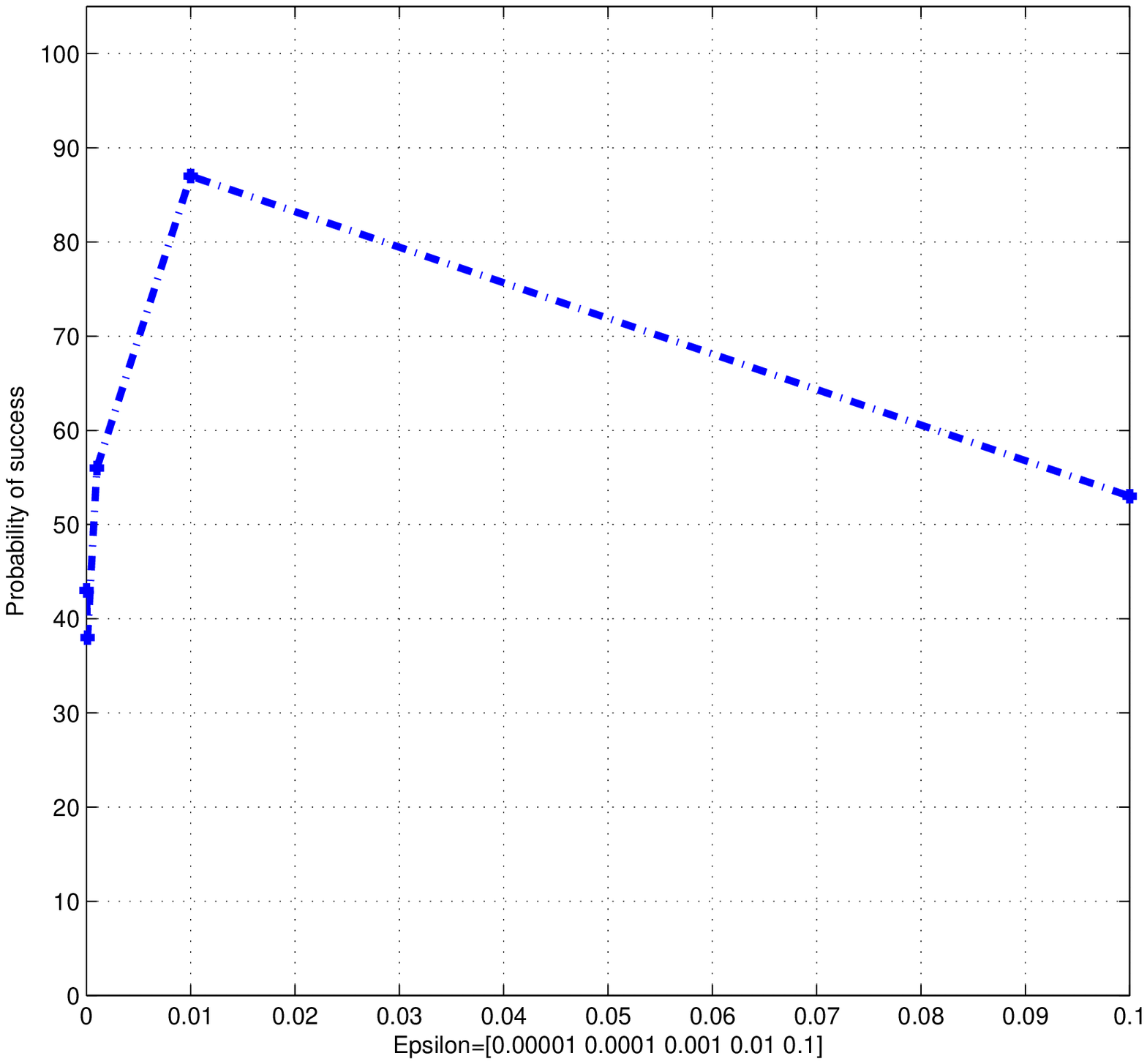}
\label{FIG14}
\end{array}$
\caption{Comparing the performance of $W_1$ minimization using different $\epsilon=0.00001,0.0001,$ $0.001,0.01,0.1$ via the probability of success for finding the exact $k=15$-sparse solution of $Ax=b$, where $A\in \mathbb{R}^{50\times 200}$, $b\in \mathbb{R}^{50}$, $p=0.05$. Matrix $A$ has been generated from Normal distribution. 100 randomly generated matrices have been tested for different chosen epsilons.}
$\begin{array}{c}
\includegraphics[width=0.5\textwidth, totalheight=.25\textheight]{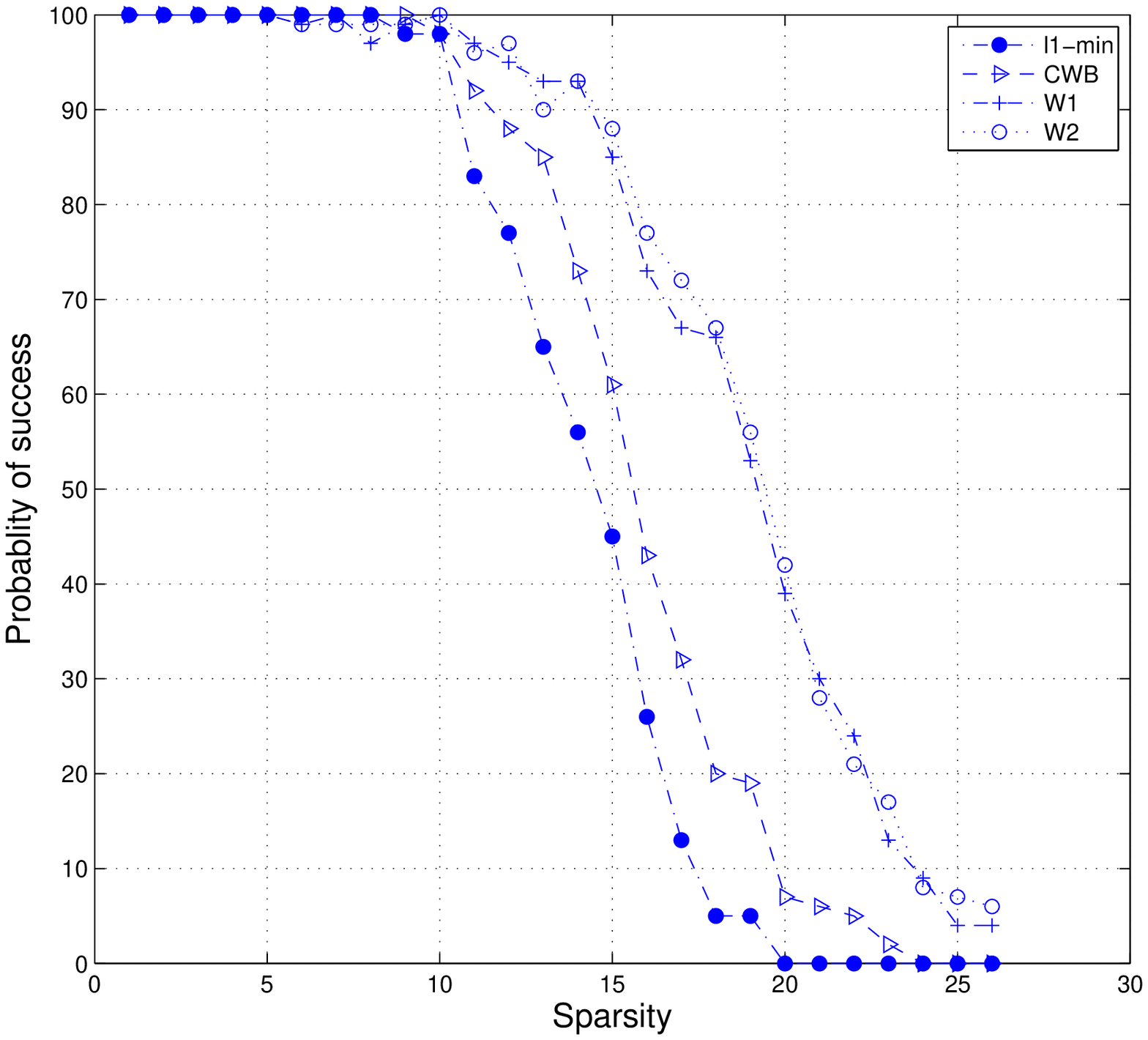}
\label{FIG14}
\end{array}$
\caption{Comparing the performance of $l_1$-min, CWB, $W_1$, $W_2$ minimization using fixed $\epsilon=0.01$ via the probability of success for finding the exact $k$-sparse solution of $Ax=b$, where $A\in \mathbb{R}^{50\times 200}$, $b\in \mathbb{R}^{50}$, $p=q=0.05$. Matrix $A$ has been generated from Normal distribution. 100 randomly generated matrices have been tested for different sparsity of $k=1,...,26$.}
\end{figure}
\begin{figure}[htp]
\centering
$\begin{array}{c}
\includegraphics[width=0.5\textwidth,totalheight=0.25\textheight]{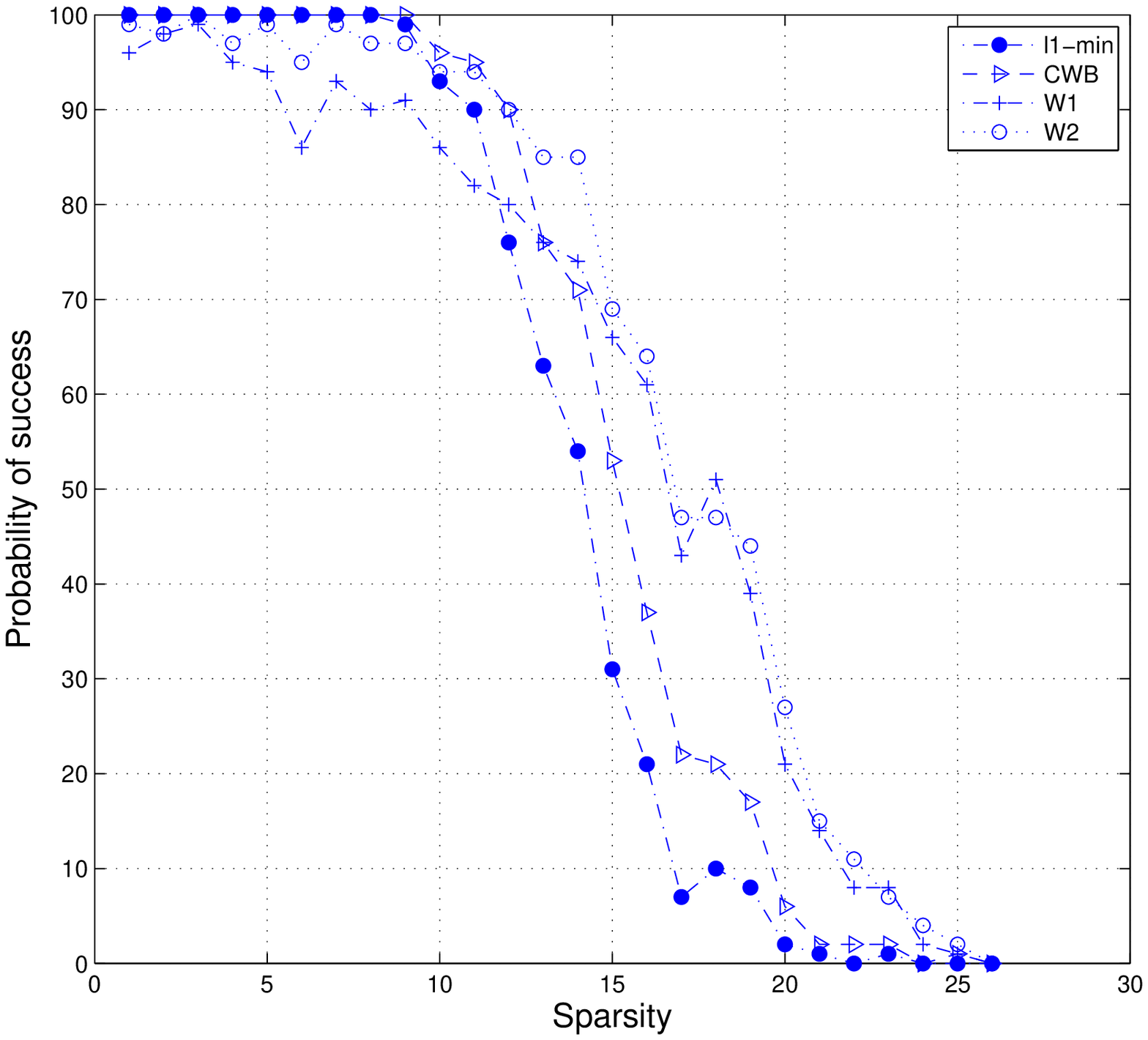}
\label{FIG13}
\end{array}$
\caption{Comparing the performance of $l_1$-min, CWB, $W_1$, $W_2$ minimization using Candes updates rule via the probability of success for finding the exact $k$-sparse solution of $Ax=b$, where $A\in \mathbb{R}^{50\times 200}$, $b\in \mathbb{R}^{50}$, $p=q=0.05$. Matrix $A$ has been generated from Normal distribution. 100 randomly generated matrices have been tested for different sparsity of $k=1,...,26$.}
\end{figure}

\bibliographystyle{amsplain}
\bibliography{Paper3}
\end{document}